\documentclass{elsart}

\usepackage{epsfig}

\usepackage{amsmath,amsfonts}

\usepackage{amssymb}

\usepackage{rotating}

\usepackage{subfigure}

\newcommand{\pd}[2]{\frac{\partial#1}{\partial#2}}
\newcommand{\pdd}[2]{\frac{{\partial#1}^2}{\partial#2}}

\begin{document}

\begin{frontmatter}

\title{A comparative study of two stochastic mode reduction methods}

\author{Panagiotis Stinis}
\address{Department of Mathematics, 
Lawrence Berkeley National Laboratory, CA 94720 USA.}
\ead{stinis@math.lbl.gov}

\begin{abstract}
We present a comparative study of two methods 
for the reduction of the 
dimensionality of a system of ordinary differential
equations that exhibits time-scale separation. Both methods lead to a reduced system of stochastic differential equations. The novel feature of these methods is that they allow the use, in the reduced system, of higher order terms in the resolved variables. The first method, proposed by Majda, Timofeyev and Vanden-Eijnden, is based on an asymptotic strategy developed by Kurtz. The second method is a short-memory approximation of the Mori-Zwanzig projection formalism of irreversible statistical mechanics, as proposed by Chorin, Hald and Kupferman. We present conditions under which the reduced models arising from the two methods should have similar predictive ability. We apply the two methods to test cases that satisfy these conditions. The form of the reduced models and the numerical simulations show that the two methods have similar predictive ability as expected.

\end{abstract}


\end{frontmatter}

\section{Introduction}
In recent years have appeared different methods for reducing the dimensionality 
of a system of ordinary or stochastic differential equations. The realization
that despite the rapid increase of computational power there are many
problems that are too expensive to tackle directly, or that in several 
problems the objects of interest are macroscopic, coarse-grained quantities, 
has led many researchers to construct methods for extracting models of
reduced dimensionality. Since the motivation to develop such methods comes
usually from specific physical problems, 
most methods exploit the mathematical structure of the problem at hand, and 
thus are not suited for all problems where dimensional reduction is needed.

Since there exists a whole range of dimension reduction methods, it is interesting 
to compare the behaviour of the reduced models predicted by different methods 
in examples that are inspired by problems of physical 
interest. The purpose of the present work is to report some results of the application 
of the asymptotic mode reduction strategy (AMRS) and of the short-memory Mori-Zwanzig 
approximation (short-memory MZ) 
to a collection of models that share common features with more complicated
systems that appear in the study of climate dynamics. These non-trivial test 
cases appeared in \cite{majda2} where an analysis of their
behaviour and of the performance of AMRS can be found. The analysis here 
suggests that the asymptotic mode reduction strategy and short-memory MZ approximation 
have, under conditions that we discuss in Section \ref{simdif}, similar predictive ability. The test cases satisfy these conditions and thus, we expect the reduced models arising from the two methods to have similar predictive ability. The form of the reduced models and the numerical simulations support our expectations.

The majority of the methods for dimensional reduction rely on a separation
of time-scales between the variables of interest and the rest of the variables 
in the system under investigation. This separation of scales can be
used in different ways depending on the form of the solutions for the fast 
dynamics. Thus, in the case of inertial manifold methods \cite{const}, the fast 
variables are, after a short transient, slaved to the slow ones. In the case
of methods that are known as "averaging methods" \cite{sand}, the fast variables
are allowed to have more complex behaviour and affect the slow variables through 
the empirical measure on the fast dynamics. The resulting equations for the slow 
variables are deterministic. In the case of the asymptotic mode reduction strategy 
 \cite{majda1}, the fast variables reach a 
statistical equilibrium much faster than the slow variables and this is used to 
obtain (in a certain limit) the reduction of the system's dimension. The resulting
system of equations for the slow variables is stochastic. In the case of 
\cite{kupf1} the resolved variables are coupled to a heat bath that is then 
approximated using the trigonometric representation for Gaussian processes. 
Chorin, Hald and Kupferman \cite{chorin1,chorin2}, 
have proposed the use of the Mori-Zwanzig projection formalism of 
irreversible statistical mechanics, for the reduction of the dimensionality of a
system of ordinary differential equations. In 
contrast to the previous methods, the Mori-Zwanzig projection formalism results  
in a reformulation of the equations for the resolved variables. This does not
make the problem of dimensional reduction easier. However, the reformulation 
serves as a starting point for approximations of various degrees of sophistication 
depending on the properties of the problem at hand. Finally, we should note the work of Papanicolaou \cite{papa,papa2,papa3,papa4} which uses a projection formalism (similar to Mori-Zwanzig, albeit with a different projection operator) on the level of the  Fokker-Planck and backward Fokker-Planck (Chapman-Kolmogorov) equations.

Every mode reduction method has two parts: a) identify the variables that will 
constitute the reduced description and b) find the equations for these variables. 
The techniques mentioned above avoid the first part by assuming that the variables to
be picked are known in advance. However, knowledge of what are the "right" variables to 
pick in an arbitrary system is a very important and difficult problem, since 
different combinations of variables can lead to very different reduced equations. 
This is an important issue for the efficient numerical implementation of any
mode reduction method. Techinques based on proper orthogonal decomposition 
\cite{siro}, coarse-grained integration \cite{kevr} and the transfer operator 
\cite{deuf} attempt not only to perform dimension reduction, but also to 
identify the appropriate variables that will constitute the reduced system. 
A concise exposition of all the methods mentioned above along with applications
to models can be found in \cite{kupf2}.

The paper is organized as follows. In Section \ref{amrs} we present a summary 
of the asymptotic mode reduction strategy as it is applied to a system of
ordinary differential equations. In Section \ref{op} we present the Mori-Zwanzig projection 
formalism and the short-memory approximation arising 
from it which is appropriate for the examples that 
we will be studying. This approximation serves to illustrate some analogies and differences between the two methods and we present them in Section \ref{simdif}. We also present conditions under which the two methods are expected to have similar behaviour. In Section \ref{models} we present the equations for the test cases that we will be examining. Section \ref{reduce} collects the results of the application of the two methods to the test cases. A short discussion follows in Section \ref{conc}.

\section{Asymptotic mode reduction strategy}{\label{amrs}}

Suppose we are given a system of ordinary differential equations
\begin{equation}
\frac{d\phi}{dt}=R(\phi) \label{odes} 
\end{equation}
with initial condition $\phi(0)=x.$ The asymptotic mode reduction strategy 
\cite{majda1,majda3} is a two step procedure based on the assumption that 
the set of variables of the system have been split into two subsets, the resolved 
and unresolved. The objective is to construct equations for the evolution of the 
resolved variables by eliminating the unresolved variables. The splitting of the variables 
in resolved and unresolved is not arbitrary but relies on the assumption that there 
exist two sets of variables in the system that evolve on different time-scales. The slow 
variables are identified as resolved ones while the fast variables as unresolved. 
In the first step, the equations 
for the unresolved variables are modified by replacing the nonlinear self-interaction 
terms between the  unresolved variables by stochastic terms. The motivation for this is 
that the self-interaction terms govern the perturbations on short time-scales, and if
we are interested in coarse-grained modelling on longer time-scales, these terms
can be replaced accurately by stochastic terms. It is important to note that this 
replacement is an assumption that has to be checked case by case. 
The second step uses 
projection techniques for stochastic differential equations \cite{kurtz1,kurtz2,ellis,papa} to eliminate the fast variables. The elimination procedure can be made rigorous 
in the limit of infinite separation of time-scales between resolved and unresolved 
variables (see e.g. \cite{kurtz1}).

We present a summary of the two steps of the method for a system of ordinary 
differential equations. We will focus on equations whose form can describe all the 
examples that we will study. 
More general forms can be found in \cite{majda1,kupf2} and references there.

For the first step we assume that we have obtained 
in some way the splitting of the variables in resolved and unresolved such that these 
two sets evolve on different time-scales. Suppose 
that the vector $\phi$ of system variables can be written as $\phi=(\hat{\phi},\tilde{\phi}),$
where $\hat{\phi}$ are the resolved variables and $\tilde{\phi}$ the unresolved ones. 
Similarly, the initial conditions can be written as $x=(\hat{x},\tilde{x}).$ The 
first step is to represent the nonlinear self-interactions between 
unresolved variables in the equations for the unresolved variables with suitable 
stochastic terms. In general, it is not easy to justify such a replacement (see e.g. \cite{just1}). 
For the examples that we study, the justification has 
been numerical \cite{majda2}. We will return to this point when 
we present the test cases. We rewrite the system (\ref{odes}) as
\begin{equation}
\label{rodes}
\begin{split}
\frac{d\hat{\phi}}{dt}&=\hat{R}(\phi) \\
\frac{d\tilde{\phi}}{dt}&=\tilde{R}({\phi})=H(\phi)+S(\tilde{\phi}),
\end{split} 
\end{equation}
where the term $S(\tilde{\phi})$ contains the nonlinear self-interaction terms 
between the unresolved variables and $H(\phi)$ contains the rest of the terms 
on the RHS of the equations for the unresolved variables. The first step of AMRS 
replaces the term $S$ by the following stochastic expression 
\begin{equation}
S(\tilde{\phi}) \approx -\frac{\Gamma}{\epsilon} \tilde{\phi}+
\frac{\sigma}{\sqrt{\epsilon}} \frac{d\vec{W}}{dt},
\label{ou}
\end{equation}
where $\Gamma,\sigma$ are positive definite (diagonal in our case) matrices and 
$\vec{W}(t)$ 
is a vector-valued Wiener process. As mentioned in \cite{majda1}, 
the choice of this particular expression in not essential to the theory but it is convenient 
for the calculations because of 
the availability of explicit expressions for the inverse of the associated backward 
Fokker-Planck operator.
The parameter 
$\epsilon$ roughly measures the ratio of the correlation times of the unresolved variables 
to the resolved ones and its precise definition will be given later when  we study 
the test cases. 
The stochastic approximation 
introduces the matrices $\Gamma,\sigma$ that have to be determined before we proceed with 
the second step of AMRS. The entries of $\Gamma$ and 
$\sigma$ are chosen so as to minimize the difference between 
the statistical moments of the original system and those of the approximate stochastic 
system (we refer to this as statistical agreement). 
For our choice of stochastic approximation which 
introduces two parameters, the procedure will involve 
the one-time and two-time statistical properties of (\ref{odes}).  We can use the 
fluctuation-dissipation theorem (one-time statistics) to reduce the number of parameters to be determined to only the entries of the matrix $\Gamma.$ The matrix $\Gamma$ is determined through the two-time statistics of (\ref{odes}). More details on how to determine $\Gamma,\sigma$ will be given when we present the test cases.

With the stochastic approximation (\ref{ou}) 
the system (\ref{rodes}) is replaced by
\begin{equation}
\label{rndodes}
\begin{split}
\frac{d\hat{\phi}}{dt}&=\hat{R}(\phi) \\
\frac{d\tilde{\phi}}{dt}&=H(\phi) -\frac{\Gamma}{\epsilon}\tilde{\phi}+
\frac{\sigma}{\sqrt{\epsilon}}\frac{d\vec{W}(t)}{dt}.
\end{split} 
\end{equation}
To facilitate the presentation, we can think of the
following splitting of the RHS for the resolved variables,
$$\hat{R}(\phi)=\hat{R}^0(\phi)+\epsilon \hat{R}^1(\hat{\phi}),$$ where
$\hat{R}^0(\phi)$ is the fast part of the RHS that includes the contributions
from the unresolved variables, while $\epsilon \hat{R}^1(\hat{\phi})$ is the
slow part that comes only from interactions between the resolved variables.
The asymptotic mode reduction strategy can be made rigorous in the limit of 
$\epsilon \rightarrow 0.$ To prepare the equations for the final elimination 
of the unresolved variables, we coarse-grain the time $t \rightarrow \epsilon t$ 
and we have
\begin{equation}
\label{cgrndodes}
\begin{split}
\frac{d\hat{\phi}}{dt}&=\frac{1}{\epsilon}\hat{R}^0(\phi)+\hat{R}^1(\hat{\phi}) \\
\frac{d\tilde{\phi}}{dt}&=\frac{1}{\epsilon}H(\phi) -
\frac{\Gamma}{{\epsilon}^2}\tilde{\phi}+
\frac{\sigma}{\epsilon}\frac{d\vec{W}(t)}{dt}.
\end{split} 
\end{equation}
As can be seen from (\ref{cgrndodes}), the dynamics for $\tilde{\phi}$ are an 
order of magnitude faster than the ones for $\hat{\phi}.$

Before we proceed with the
second step of AMRS that eliminates the unresolved variables, 
we note again that the approximation (\ref{ou}) is a 
working assumption that has to be checked case by case. However, if we assume that 
the approximation (\ref{ou}) is valid, then, due to a theorem by Kurtz \cite{kurtz1}, 
the process of elimination of the unresolved variables becomes rigorous in the 
limit of inifinite scale separation. 
The theorem operates on the level of the Chapman-Kolmogorov (or backward Fokker-Planck) 
equation associated with (\ref{cgrndodes}).

When presented with a system of stochastic differential equations of the form 
(\ref{cgrndodes}) we can construct two linear partial differential equations, 
one for the evolution of the probability density function 
(Fokker-Planck equation) and one for the evolution of expectation values of 
functions of the solution (Chapman-Kolmogorov equation). The 
Chapman-Kolmogorov equation is the adjoint of the Fokker-Planck equation 
with respect to the scalar product defined on  $X=\mathbb{R}^n,$ where $n$ is the 
dimensionality of the system (\ref{odes}) (see below). 
We shall focus our attention on the 
Chapman-Kolmogorov equation, since we are only interested in giving a brief account 
of AMRS (see \cite{kupf2} for more details). For any scalar-valued function 
$g(\hat{x})$, we define the quantity 
$$u^{\epsilon}(x,t)=\mathbf{E} g(\hat{\phi}).$$ Then 
$u^{\epsilon}$ satisfies the Chapman-Kolmogorov partial differential equation 
\begin{equation}
\label{ck}
\frac{\partial{u^{\epsilon}}}{\partial{t}}=
\frac{1}{{\epsilon}^2}\mathcal{L}_1 u^{\epsilon}
+\frac{1}{\epsilon}\mathcal{L}_2 u^{\epsilon} 
+\mathcal{L}_3u^{\epsilon},
\end{equation}
with $u^{\epsilon}(x,0)=g(\hat{x}).$ Note that while $u^{\epsilon}(x,0)$ is only a 
function of $\hat{x}$, $u^{\epsilon}(x,t)$ is a function of all the variables. 
The operators $\mathcal{L}_1,\mathcal{L}_2$ and $\mathcal{L}_3$ 
are defined by
\begin{equation}
\label{ckoper}
\begin{split}
\mathcal{L}_1&=\sum_{j \in I_{\tilde{x}}} 
\biggl( -\gamma_j x_j \pd{}{x_j} + \frac{{\sigma_j}^2}{2} \pdd{}{x_j^2}\biggr)\\
\mathcal{L}_2&=\sum_{j \in I_{\hat{x}}} \hat{R}^0_j(x) \pd{}{x_j}
+\sum_{j \in I_{\tilde{x}}} H_j(x) \pd{}{x_j} \\
\mathcal{L}_3&=\sum_{j \in I_{\hat{x}}} \hat{R}^1_j(\hat{x}) \pd{}{x_j}
\end{split}
\end{equation}
where $I_{\hat{x}},I_{\tilde{x}}$ are index sets for the resolved and  
unresolved variables respectively. Also, $\gamma_j, \sigma_j$ are the (diagonal) 
entries of the matrices $\Gamma, \sigma.$

The second step of AMRS eliminates the unresolved variables. From (\ref{ck}) we 
see that the limit $\epsilon \rightarrow 0$ is singular. Also, that (\ref{ck}) 
is a partial differential equation involving both resolved and unresolved 
variables. Thus, our aim is to derive (if possible), 
in the limit $\epsilon \rightarrow 0$, a Chapman-Kolmogorov equation for the resolved 
variables only and from this equation read off the form of the corresponding 
stochastic system for the resolved variables. We 
will use only formal manipulations that can be made rigorous in the limit 
of $\epsilon \rightarrow 0$ through Kurtz's theorem. We proceed by 
expanding the solution $u^{\epsilon}(x,t)$ in powers of $\epsilon$
\begin{equation}
\label{cksper}
u^{\epsilon}(x,t)=u_0+\epsilon u_1+ {\epsilon}^2 u_2 + \ldots.
\end{equation}
We insert (\ref{cksper}) into (\ref{ck}) and collect terms of equal power in $\epsilon.$ 
We obtain the following hierarchy of equations
\begin{align}
\mathcal{L}_1 u_0&=0 \label{cksper1}\\
\mathcal{L}_1 u_1&=-\mathcal{L}_2 u_0  \label{cksper2}\\
\mathcal{L}_1 u_2&=\pd{u_0}{t}-\mathcal{L}_2 u_1 \label{cksper3}-\mathcal{L}_3 u_0 \\
         & \vdots \notag
\end{align}
Define a scalar product on $X$ by
\begin{equation}
\label{scalar}
\langle a,b \rangle=\int_{X}a(x)b(x)dx,
\end{equation}
for any two scalar-valued functions $a,b.$
The solvability of the equations (\ref{cksper1})-(\ref{cksper3}) requires that the 
RHS of the equations belong to the range of $\mathcal{L}_1$, or equivalently that 
they are orthogonal (under the scalar product just defined) 
to the kernel of $\mathcal{L}_1^{*},$ the adjoint of 
$\mathcal{L}_1.$ By construction, we have chosen the dynamics of $\mathcal{L}_1$ in 
such a way that the kernel of $\mathcal{L}_1^{*}$ contains only one element. This is the 
invariant density $f_s(\tilde{x})$ under the dynamics governed by the linear stochastic operator  in (\ref{ou}). Thus, the solvability condition for all the equations 
implies that the RHSs average to zero with respect to $f_s(\tilde{x}).$ 
Consequently, if we denote by $\mathbf{P}$ the projection on the subspace 
$\hat{X}$ of resolved variables, where 
$$\mathbf{P} \cdot \equiv \int_{\tilde{X}} \cdot f_s(\tilde{x}) d\tilde{x},$$
we have $\mathbf{P}(RHS)=0.$ Also, the solvability of the 
equation (\ref{cksper1}) implies $\mathbf{P}u_0=u_0,$ i.e. $u_0$ is independent of $\tilde{x}$ 
(see \cite{kupf2}). From equation (\ref{ckoper}) and the definition of the projection 
 $\mathbf{P}$ we have that $\mathbf{P}\mathcal{L}_1=\mathcal{L}_1\mathbf{P}=0$ 
(for more details see \cite{majda1,kupf2}).
We collect these remarks and find that the function $u_0$ satisfies the equation
\begin{equation}
\label{ckf}
\begin{split}
\pd{u_0}{t}&= \mathbf{P}\mathcal{L}_3\mathbf{P}u_0
-\mathbf{P}\mathcal{L}_2\mathcal{L}_1^{-1}\mathcal{L}_2\mathbf{P}{u_0} \\
u_0(\hat{x},0)&=g(\hat{x})
\end{split}
\end{equation}
where $\mathcal{L}_1^{-1}$ is the inverse of the operator $\mathcal{L}_1$ (see Appendix 
A in \cite{majda1}). In the limit of $\epsilon \rightarrow 0,$ Kurtz's theorem makes these 
manipulations rigorous by stating that under the condition 
$\mathbf{P}\mathcal{L}_2\mathbf{P}=0,$ i.e the solvability condition for (\ref{cksper2}), 
the solution $u^{\epsilon}(x,t)$ tends to $u_0(\hat{x},t)$, where $u_0$ satisfies (\ref{ckf}). 
Moreover, equation (\ref{ckf}) is a Chapman-Kolmogorov equation involving only the 
resolved variables and since $\mathcal{L}_1^{-1}$ is computable for the stochastic 
approximation (\ref{ou}), we can use (\ref{ckf}) to 
compute the corresponding stochastic differential equations for the resolved variables. This 
concludes the second step of AMRS. Note, that in general, the effective compuational use of the Chapman-Kolmogorov 
equation for the resolved variables relies on the operator $\mathcal{L}_1$ having a closed 
form expression for its inverse.

\section{Mori-Zwanzig projection formalism and the short-memory approximation}
{\label{op}}

We present the Mori-Zwanzig formalism and
derive from it the approximation that we will be using later, namely
the short-memory approximation. 
\subsection{Conditional expectations and the Mori-Zwanzig formalism}
Suppose that we are given the system of equations (\ref{odes}) with initial
condition $\phi(0)=x.$ Furthermore, assume that we know only
a fraction of the initial data, say $\hat{x}$,
where $x=(\hat{x},\tilde{x})$ and correspondingly
$\phi=(\hat{\phi},\tilde{\phi})$ and that the
unresolved data are drawn from a measure
with density $f(x).$ Unlike the case of AMRS, here we do not assume 
anything about the form of the density $f(x).$

Suppose $u,v$ are functions of $x$,
and introduce the scalar product $(u,v)=E[uv]=\int{u(x)}v(x)f(x)dx$.
We will denote the space of functions $u$ with $E[u^2]< \infty$ 
by $L^2(f)$ or simply $L^2$.
We are looking for approximations of functions of $x$ by functions of
$\hat{x}$, where $\hat{x}$ are the variables that form our reduced system
(the resolved degrees of freedom). The functions of $\hat{x}$ form a
closed linear subspace of $L^2$, which we denote by $\hat{L}^2$.
Given a function $u$ in $L^2$, its conditional expectation with respect to
$\hat{x}$ is
$$E[u|\hat{x}]=\frac{\int{ufd\tilde{x}}}{\int fd\tilde{x}}.$$
The conditional expectation $E[u|\hat{x}]$ is the best approximation of 
$u$ by a function of $\hat{x}$:
$$E[|u-E[u|\hat{x}]|^2] \leq E[|u-h(\hat{x})|^2]$$
for all functions $h$.

We pick a basis in $\hat{L}^2$,
for example $h_1(\hat{x}),h_2(\hat{x}),\ldots$. For simplicity assume
that the basis functions $h_i(\hat{x})$ are orthonormal, i.e.,
$E[h_ih_j]=\delta_{ij}$. The conditional
expectation can be written as
$E[u|\hat{x}]=\sum{a}_j h_j(\hat{x})$, where $a_j=E[u h_j]=
E[u(\hat{x},\tilde{x})h_j(\hat{x})]$. If we have a finite number
of terms only, we are projecting on a smaller subspace and the
projection is called a finite-rank projection. In the special case
where we pick a the finite set of functions
$h_1(\hat{x})=x_1,h_2(\hat{x})=x_2,\ldots,h_m(\hat{x})=x_m$,
then the corresponding finite-rank projection is called in physics
the "linear" projection (note that all projections are linear, so "linear" is used to
denote that the projection is on linear functions of the resolved variables). 
We should note here that it is not always true that $E[x_i x_j]=0$ for $i\neq{j}$.

The system of ordinary differential equations
we are asked to solve can be transformed into the linear
partial differential equation \cite{chorin4}
\begin{equation}
\label{pde}
u_t=Lu, \qquad u(x,0)=g(x)
\end{equation}
where $L=\sum_i R_i(x)\frac{\partial}{\partial{x_i}}$ and the solution of (\ref{pde}) is
given by $u(x,t)=g(\phi(x,t))$. Consider the following
initial condition for the PDE
$$g(x)=x_j \Rightarrow  u(x,t)=\phi_j(x,t)$$
Using semigroup notation we can rewrite (\ref{pde}) as
$$\frac{\partial}{\partial{t}} e^{tL}x_j=L e^{tL}x_j$$
Let $Q=I-P.$ Equation (\ref{pde}) 
can be rewritten as \cite{chorin4}
\begin{equation}
\label{mz}
\frac{\partial}{\partial{t}} e^{tL}x_j=
e^{tL}PLx_j+e^{tQL}QLx_j+
\int_0^t e^{(t-s)L}PLe^{sQL}QLx_jds,
\end{equation}
where we have used Dyson's formula
\begin{equation}
\label{dyson1}
e^{tL}=e^{tQL}+\int_0^t e^{(t-s)L}PLe^{sQL}ds.
\end{equation}
Equation (\ref{mz}) is the Mori-Zwanzig identity \cite{zwan1,zwan2,mori1}. 
Note that
this relation is exact and is an alternative way
of writing the original PDE. It is the starting
point of our approximations. Of course, we
have one such equation for each of the resolved
variables $\phi_j, j=1,\ldots,m$. The first term in (\ref{mz}) is
usually called Markovian since it depends only on the values of the variables
at the current instant, the second is called "noise" and the third "memory". 
The meaning of the different terms appearing in (\ref{mz}) and a connection 
(and generalization) to the fluctuation-dissipation theorems of irreversible 
statistical mechanics can be found in \cite{chorin1,stinis}. 

If we write
$$e^{tQL}QLx_j=w_j,$$ 
$w_j(x,t)$ satisfies the equation
\begin{equation}
\label{ortho}
\begin{cases}
&\frac{\partial}{\partial{t}}w_j(x,t)=QLw_j(x,t) \\ 
& w_j(x,0) = QLx_j=R_j(x)-PR_j(\hat{x}). 
\end{cases} 
\end{equation}
If we project (\ref{ortho})
using any of the projections discussed we get
$$P\frac{\partial}{\partial{t}}w_j(x,t)=
PQLw_j(x,t)=0,$$
since $PQ=0$. Also for the initial condition
$$Pw_j(x,0)=PQLx_j=0$$
by the same argument. Thus, the solution
of (\ref{ortho}) is at all times orthogonal
to the range of $P.$ We call
(\ref{ortho}) the orthogonal dynamics equation. Since the solutions of 
the orthogonal dynamics equation remain orthogonal to the subspace $\hat{L}^2$, 
we can project the Mori-Zwanzig equation (\ref{mz}) on $\hat{L}^2$ and find
\begin{equation}
\label{mzp}
\frac{\partial}{\partial{t}} Pe^{tL}x_j=
Pe^{tL}PLx_j+
P\int_0^t e^{(t-s)L}PLe^{sQL}QLx_jds,
\end{equation}

\subsection{The short-time and short-memory approximations}{\label{smop}}
The approximation we will 
examine is a short-time approximation and consists of dropping 
the integral term in Dyson's formula (\ref{dyson1}) 

\begin{equation}
e^{tQL} \approx e^{tL}  \label{sm1}.
\end{equation}
In other words we replace the flow in the orthogonal complement 
of $\hat{L}^2$ with the flow induced by the full system operator 
$L$. Some algebra shows that 
\begin{equation}
Q(e^{sQL}-e^{sL})=O(s^2) \label{sm41},
\end{equation}
and
\begin{equation}
\int_0^t e^{(t-s)L}PLe^{sQL}QLx_jds=\int_0^t e^{(t-s)L}
PLQe^{sL}QLx_jds+O(t^3). \label{sm5}
\end{equation}

As expected the approximation (\ref{sm1}) is good only for short times. However, 
under certain conditions this approximation can become valid 
for longer times. To see that consider the case where $P$ is 
the finite-rank projection so
$PLQe^{sQL}QLx_j =\sum_{k=1}^l (LQe^{sQL}QLx_j,h_k)h_k(\hat{x})$
and $PLQe^{sL}QLx_j =\sum_{k=1}^l (LQe^{sL}QLx_j,h_k)h_k(\hat{x}).$
The quantities $(LQe^{sL}QLx_j,h_k)$ can be calculated from 
the full system without recourse to the orthogonal dynamics. 
Recall (\ref{sm41}) which states that the error in 
approximating $e^{sQL}$ by $e^{sL}$ is small for small $s$. 
This means that for short times we can infer the behaviour of 
the quantity $(LQe^{sQL}QLx_j,h_k)$ by examining 
the behaviour of  the quantity $(LQe^{sL}QLx_j,h_k)$.

If the quantities $(LQe^{sL}QLx_j,h_k)$ decay fast we can infer 
that the quantities $(LQe^{sQL}QLx_j,h_k)$ decay 
fast for short times. We cannot infer anything about the
behaviour of $(LQe^{sQL}QLx_j,h_k)$ for
larger times. However, if $(LQe^{sQL}QLx_j,h_k)$ not
only decay fast initially, but, also, stay small for larger times,
then we expect our approximation to be valid for larger
times. To see this consider again the integral term in the 
Mori-Zwanzig equation. The integrand becomes negligible for $ t \gg t_0$, 
where $t_0$ is the time of decay of the quantities 
$(LQe^{sQL}QLx_j,h_k)$. 
This means that our approximation becomes
$$\int_0^t e^{(t-s)L}PLe^{sQL}QLx_jds \approx 
\int_{0}^{t_0} e^{(t-s)L}PLQe^{sL}QLx_jds+O(t_0^3) \label{sm8}.$$
From this we conclude that the short-time approximation is 
valid for large times if $t_0$ is small and is called the 
short-memory approximation. On the other hand, if $t_0$ is 
large, then the error, which is $O(t_0^3)$, becomes large and 
the approximation is only valid for short times. Note that the
validity of the short-memory approximation can only be
checked after constructing it, since it
is based on an assumption about the large time behaviour
of the unknown quantities $(LQe^{sQL}QLx_j,h_k)$. Note, that
determination of the quantities $(LQe^{sQL}QLx_j,h_k)$ requires 
the (usually very expensive) solution of the orthogonal dynamics equation.
The short-memory approximation, when valid, allows us to avoid the
solution of the orthogonal dynamics equation.

If the quantities $(LQe^{sL}QLx_j,h_k)$ do not decay
fast, then we can infer, again only for short times, that the quantities 
$(LQe^{sQL}QLx_j,h_k)$ of the exact Mori-Zwanzig equation do not
decay fast. Yet, it is possible that the quantities 
$(LQe^{sQL}QLx_j,h_k)$  start decaying very fast after short times and remain
small for longer times, 
so that the short-memory approximation could still hold. Of course, 
this can only be checked a posteriori, after the simulation of the short-memory
approximation equations.

In the statistical physics literature, the assumption that the 
correlations vanish for $s\neq 0$ is often made which is a special 
case of the short-memory approximation with the correlations replaced 
by a delta-function multiplied by the integrals. An application of the 
short-memory approximation can be found in \cite{bell1}.

\section{A formal comparison of the two methods}{\label{simdif}}

We use the conditions under which the asymptotic mode reduction 
strategy holds to derive conditions under which we can construct a short-memory 
reduced model in the context of the MZ formalism. This allows us to obtain 
a form of the reduced system in the MZ formalism that is more amenable to comparison 
to the reduced model obtained by AMRS. We should note here that it is possible to start from the Chapman-Kolmogorov equation (\ref{ck}) and apply a projection formalism similar to Mori-Zwanzig using the operator $\mathbf{P}$ not $P$ (see e.g. \cite{papa2} and references therein). Under suitable assumptions (the same conditions as the ones we will postulate for the derivation of the short-memory MZ model in this section) one can arrive at the reduced Chapman-Kolmogorov equation (\ref{ckf}). Our intention in this section is different. We want to compare the reduced short-memory model that arises when MZ is applied to the {\it deterministic} system (\ref{odes}) to the reduced short-memory model that arises from AMRS, where the system (\ref{odes}) is first replaced by the {\it stochastic} system (\ref{rodes}).

Write $S=\frac{1}{\epsilon}S_0.$ This is along the lines of introducing explicitly the scaling of the stochastic operator in (\ref{ou}) (note that in the following we continue to refer to the dynamics of the term as $S$ dynamics and not $S_0$ dynamics). If, in addition, we coarse-grain time $t \rightarrow \epsilon t,$ the system 
(\ref{odes}) becomes 
\begin{equation}
\label{odesr}
\begin{split}
\frac{d\hat{\phi}}{dt}&=\frac{1}{\epsilon}\hat{R}^0(\phi)+\hat{R}^1(\hat{\phi}) \\
\frac{d\tilde{\phi}}{dt}&=\frac{1}{\epsilon}H(\phi) + \frac{1}{{\epsilon}^2} S_0(\tilde{\phi}).
\end{split} 
\end{equation}
The operator $L$ associated with (\ref{odesr}) is 
\begin{equation}
\label{odesrl}
L=\frac{1}{{\epsilon}^2}L_1
+\frac{1}{\epsilon}L_2 
+L_3
\end{equation}
where 
\[
L_1=\sum_{j \in I_{\tilde{x}}} S_{0_j}(x) \pd{}{x_j}, \quad
 L_2=\mathcal{L}_2 \quad \text{and} \quad L_3=\mathcal{L}_3.
\]
From (\ref{mzp}) we find for the resolved variable $x_j, j=1,\ldots,m$ 
\begin{equation}
\label{mzps}
\frac{\partial}{\partial{t}} Pe^{tL}x_j=
Pe^{tL}PLx_j+
P\int_0^t e^{(t-s)L}PLQe^{sQL}QLx_jds,
\end{equation}
where $L$ is given by (\ref{odesrl}). We want to use the conditions under which AMRS holds 
to simplify the equation (\ref{mzps}). Note that AMRS works with the projection operator $\mathbf{P}$ while the MZ formalism works with the projection operator $P$ (we will say more about the relation between the two operators below). The asymptotic mode reduction strategy  has two steps. In the first step, under suitable assumptions, the dynamics described by $S$ are replaced by a linear stochastic operator. In the second step, under suitable assumptions which guarantee the applicability of Kurtz's theorem, the unresolved variables are eliminated and a reduced model for the resolved variables is obtained. Our goal here is to translate all these assumptions (conditions) in the language of the MZ formalism, and, thus, obtain the conditions under which, for a system with $L$ given by (\ref{odesrl}), it is possible to construct a short-memory model  for the slow variables.

The conditions under which Kurtz's theorem holds are, as already mentioned in Section \ref{amrs}, 
\begin{equation}
\label{amrscon}
\mathbf{P}\mathcal{L}_2\mathbf{P}=0, \quad \mathbf{P}\mathcal{L}_1=\mathcal{L}_1\mathbf{P}=0.
\end{equation}
For the operator $P$ used in the MZ formalism we have , by construction, $L_1P=0.$ In the same way that AMRS can be applied if the conditions (\ref{amrscon}) hold, we suggest that similar conditions should hold for the operator $P$, if a short-memory model is to be found for the special case where $L$ is given by (\ref{odesrl}). Thus, we begin by postulating the conditions
\begin{equation}
\label{mzcon}
PL_2P=0, \quad PL_1=L_1P=0.
\end{equation}
The condition $PL_2P=0$ implies that $Pe^{tL}PL_2Px_j=0.$ We note here that the 
conditions (\ref{mzcon}) hold for the test cases. Under these conditions, the equation for the resolved variable $x_j$ becomes
\begin{equation}
\label{mzps1}
\frac{\partial}{\partial{t}} Pe^{tL}x_j=
Pe^{tL}PL_3Px_j+
\frac{1}{\epsilon^2}P\int_0^t e^{(t-s)L}PL_2e^{sQL}L_2Px_jds,
\end{equation}
where we have also used the fact that $Px_j=x_j.$ The first term is the analog of the first term in 
(\ref{ckf}) and is the contribution of the interactions among the resolved variables to the reduced model. The second term in (\ref{mzps1}) comes from the interaction among the resolved and 
unresolved variables and is the analog of the second term in (\ref{ckf}). The analogy becomes more clear when we examine the meaning of the integrand  $PL_2e^{sQL}L_2Px_j.$ The existence (or not) of a reduced model in the limit $\epsilon \rightarrow 0$ depends on the behaviour of the integrand $PL_2e^{sQL}L_2Px_j$ and, in particular, on the behaviour of $L_2e^{sQL}L_2Px_j$. To see this better we can consider the case when $P$ is the finite-rank projection on a set of orthonormal functions $h_1(\hat{x}), 
\ldots,h_l(\hat{x}).$ The integral term in (\ref{mzps1}) becomes
\[
\frac{1}{\epsilon^2}P\int_0^t\sum_{k=1}^l (L_2e^{sQ(\frac{1}{{\epsilon}^2}L_1
+\frac{1}{\epsilon}L_2 
+L_3)}L_2x_j,h_k)h_k(\hat{\phi}(t-s))ds
\]
The kernels $(L_2e^{sQ(\frac{1}{{\epsilon}^2}L_1+\frac{1}{\epsilon}L_2 +L_3)}L_2x_j,h_k)$ are 
correlations of the resolved-unresolved interaction terms under the orthogonal dynamics. Given the specific structure of the operator $L$ and the conditions (\ref{mzcon}), the kernels $(L_2e^{sQ(\frac{1}{{\epsilon}^2}L_1+\frac{1}{\epsilon}L_2 +L_3)}L_2x_j,h_k)$ can be approximated, in the limit 
$\epsilon \rightarrow 0$, by $(L_2e^{s\frac{1}{{\epsilon}^2}L_1}L_2x_j,h_k).$ So, in the limit of $\epsilon \rightarrow 0,$ the orthogonal dynamics are dominated by the $S$ dynamics of the unresolved variables. This is also the meaning of the operator $\mathcal{L}_2\mathcal{L}_1^{-1}\mathcal{L}_2$ in  (\ref{ckf}). 
 Moreover, if the correlations  $(L_2e^{s\frac{1}{{\epsilon}^2}L_1}L_2x_j,h_k)$ decay on a time-scale $O(\epsilon^2)$, then the integral term contributes to the RHS of (\ref{mzps1}) at the same order with the first term. Equation (\ref{mzps1}) becomes 
\begin{equation}
\label{mzps2}
\frac{\partial}{\partial{t}} Pe^{tL}x_j=
Pe^{tL}PL_3Px_j+
\frac{1}{\epsilon^2}P\int_0^t e^{(t-s)L}PL_2e^{s\frac{1}{{\epsilon}^2}L_1}L_2Px_jds,
\end{equation}  
Equation (\ref{mzps2}) is, for the specific $L$ given by (\ref{odesrl}), the form that the short-memory MZ approximation equation for the conditional expectation of the resolved variable $x_j$ acquires.

In summary, for $L$ given by (\ref{odesrl}), we have that: i) the conditions given by (\ref{mzcon}), ii) the approximation $(L_2e^{sQ(\frac{1}{{\epsilon}^2}L_1+\frac{1}{\epsilon}L_2 +L_3)}L_2x_j,h_k) \approx  (L_2e^{s\frac{1}{{\epsilon}^2}L_1}L_2x_j,h_k)$ and iii) the assumption that the quantities $(L_2e^{s\frac{1}{{\epsilon}^2}L_1}L_2x_j,h_k)$ decay fast, guarantee the existence of a short-memory MZ model in the limit $\epsilon \rightarrow 0.$ The reduced model (\ref{mzps2}) bears great formal similarity to the reduced model (\ref{ckf}) obtained by AMRS (see also Appendix A in \cite{majda1} for the integral representation of  $\mathcal{L}^{-1}$). However, there are two main points of difference that we now discuss.

The first point is the difference between the operators $L_1$ and $\mathcal{L}_1.$ This difference arises from the stochastic approximation (\ref{ou}) in the first step of AMRS. As we have already mentioned, it is not easy to justify this approximation in the general case \cite{just1}.  For the test cases the justification was numerical \cite{majda2}.

The second point is the difference between the projection operators $P$ and $\mathbf{P}.$ The operator $P$ is defined by 
$$P \cdot \equiv \int_{\tilde{X}}  \cdot f_c d\tilde{x},$$
where $f_c(\hat{x},\tilde{x})=\frac{f}{\int_{\tilde{X}} fd\tilde{x}}$ is the conditional density conditioned on the resolved variables $\hat{x}.$ On the other hand, the operator $\mathbf{P}$ is defined by 
$$\mathbf{P} \cdot \equiv \int_{\tilde{X}} \cdot f_s d\tilde{x},$$
where $f_s(\tilde{x})$ is the invariant density of the dynamics governed by the linear stochastic operator appearing in (\ref{ou}). Thus, the question if and how $P$ relates to $\mathbf{P}$ is reduced to the question of if and how $f_c$ is related to $f_s.$ It is not obvious, if and how $f_c$ and $f_s$ are related in general. However, for the case when the original system (\ref{odes}) admits an invariant density $f,$ the question is ultimately connected to the first point mentioned above, i.e. whether the replacement of the $S$ dynamics by a stochastic operator is justified. To see that we will need a result for the conditional density $f_c$ in the limit $\epsilon \rightarrow 0,$ which we now derive. The result states that, if the density $f$ is invariant for the system (\ref{odes}), then, in the limit $\epsilon \rightarrow 0,$ the conditional density $f_c$ is invariant under the $S$ dynamics.

Indeed, if the density $f$ is invariant for the original system (\ref{odes}), then the equation for the evolution of $f$ is
\begin{equation}
\label{liouv}
L^{*}f=[\frac{1}{{\epsilon}^2}L^{*}_1+\frac{1}{\epsilon}L^{*}_2 +L^{*}_3]f=0,
\end{equation}  
where $L^{*}$ is the adjoint of $L$ with respect to the scalar product $\langle,\rangle$ defined in (\ref{scalar}). Note that, by construction 
$$ L^{*}_2=\mathcal{L}^{*}_2 \quad \text{and} \quad L^{*}_3=\mathcal{L}^{*}_3$$
but, in general, 
$$L^{*}_1 \neq \mathcal{L}^{*}_1.$$
In the limit $\epsilon \rightarrow 0$ and, if we assume that the measure remains smooth enough to admit a density, the term $\frac{1}{{\epsilon}^2}L^{*}_1f$ should dominate and (\ref{liouv}) becomes
\begin{equation}
\label{liouv2}
L_1^{*}f=0
\end{equation}
By the definition of $f_c$ we have $f=f_c \int_{\tilde{X}} fd\tilde{x}$ and (\ref{liouv2}) gives
\begin{equation}
\label{liouv3}
L_1^{*}[f_c  \int_{\tilde{X}} fd\tilde{x}]= [\int_{\tilde{X}} fd\tilde{x}]L_1^{*}f_c=0,
\end{equation}
since $L^{*}_1$ contains derivatives only with respect to the unresolved variables $\tilde{x}$ while 
$\int_{\tilde{X}} fd\tilde{x}$ is a function of the resolved variables $\hat{x}.$ From (\ref{liouv3}) we find 
\begin{equation}
\label{liouv4}
L_1^{*}f_c=0.
\end{equation}
Equation (\ref{liouv4}) is the result we need. It states that, in the limit $\epsilon \rightarrow 0$, the conditional density $f_c$ is an invariant density for the $S$ dynamics. On the other hand, we have, by definition,  
\begin{equation}
\label{liouv5}
\mathcal{L}^{*}_1f_s=0.
\end{equation}
If we assume that the equations (\ref{liouv4}-\ref{liouv5}) admit each a unique solution (for the same boundary conditions), then the relation between the densities $f_c$ and $f_s$ is governed by the relation between the operators $L^{*}_1$ and $\mathcal{L}^{*}_1.$ But the relation between   $L^{*}_1$ and $\mathcal{L}^{*}_1$ is determined by the relation between the operators $L_1$ and $\mathcal{L}_1,$ i.e. by whether or not we are allowed to replace the fast dynamics $S$ by a stochastic operator. So, the relation between the two projection operators $P$ and $\mathbf{P}$, which are defined by $f_c$ and $f_s$ respectively, is determined, in the case that $f$ is invariant for (\ref{odes}), by whether or not we are allowed to replace the fast dynamics $S$ by a stochastic operator.  For the test cases we have $f_s=f_c$ and thus, we expect the reduced models constructed from the two methods to have similar predictive ability. The numerical simulations support our expectations.

\section{The models}{\label{models}}

The models that we use to compare the two stochastic mode reduction methods mentioned 
above first appeared in \cite{majda2}, where the performance of the asymptotic mode 
reduction strategy was evaluated. We chose these examples since, there exist 
already published results about them for the asymptotic mode reduction strategy and 
all the parameters involved are well-documented. All the models that we shall examine 
have a common structure. They consist of a system of deterministic differential equations 
where a few (two at most) slowly evolving variables are coupled to a fast evolving heat 
bath. The heat bath comes from a Fourier-Galerkin truncation of the Burgers-Hopf system 
$$u_t+\frac{1}{2}(u^2)_x=0,$$
where $x \in [0,2\pi].$ 
If we expand the solution $u(x,t)$ in Fourier series 
$$u(x,t)=\sum_{k=1}^{\Lambda} u_k(t) e^{ikx},$$
we find for the mode $u_k$   
$$\frac{du_k}{dt}=-\frac{ik}{2}\sum_{|k'| \leq \Lambda} u_{k'} u_{k-k'}, \quad  
1\leq k \leq \Lambda,$$
where $\Lambda$ controls the size of the truncation. Also, $u_{-k}=u_k^{*}$ and 
we set in all calculations $u_0=0.$  

As demonstrated numerically by Majda and Timofeyev in \cite{majda4}, this 
truncation (for a large enough number of modes) is a deterministic but chaotic and mixing 
system, that is ergodic on suitably defined equi-energy surfaces, and the time-correlations 
of the modes obey a simple scaling law. These properties qualify this system as a 
good candidate for a stochastic heat bath. The modes of the Galerkin truncation will be 
the unresolved variables that we wish to eliminate, thus obtaining a stochastic system 
for the slowly evolving resolved variables. The coupling of the variables on the 
RHS of the equations will be of triad type, i.e. no variable appears on the RHS 
of its own equation.

The first model is called the additive case and is given by
\begin{equation}
\label{add}
\begin{split}
\frac{dx_1}{dt}&=\lambda\sum_k b_k^{1|yz}y_k z_k, \quad 
\frac{dy_k}{dt}=-Re\frac{ik}{2}\sum_{|k'| \leq \Lambda} u_{k'} u_{k-k'}+
\lambda b_k^{y|1z}x_1 z_k, \\
\frac{dz_k}{dt}&=-Im\frac{ik}{2}\sum_{|k'| \leq \Lambda} u_{k'} u_{k-k'}+
\lambda b_k^{z|1y}x_1 y_k,
\end{split}
\end{equation}
where $u_k=y_k + i z_k.$ The interaction coefficients are of order 1 and 
satisfy
\begin{equation}
\label{addc}
b_k^{1|yz}+b_k^{y|1z}+b_k^{z|1y}=0.
\end{equation}
The parameter $\lambda$ is a measure of the strength of the coupling between the 
variable $x_1$ and the Burgers bath. The characterization additive is given in 
anticipation of the form of the reduced model for $x_1$ which will have a noise term 
of additive type.

The second model is called the multiplicative case and is given by
\begin{equation}
\label{mul}
\begin{split}
\frac{dx_1}{dt}&=\lambda\sum_k (b_k^{1|2y}x_2 y_k+b_k^{1|2z}x_2 z_k), \quad 
\frac{dx_2}{dt}=\lambda\sum_k (b_k^{2|1y}x_1 y_k+b_k^{2|1z}x_1 z_k), \\
\frac{dy_k}{dt}&=-Re\frac{ik}{2}\sum_{|k'| \leq \Lambda} u_{k'} u_{k-k'}+
\lambda b_k^{y|12}x_1 x_2, \\ 
\frac{dz_k}{dt}&=-Im\frac{ik}{2}\sum_{|k'| \leq \Lambda} u_{k'} u_{k-k'}+
\lambda b_k^{z|12}x_1 x_2,
\end{split}
\end{equation}
where the interaction coefficients are of order 1 and 
satisfy
\begin{equation}
\label{mulc}
b_k^{1|2y}+b_k^{2|1y}+b_k^{y|12}=0, \quad b_k^{1|2z}+b_k^{2|1z}+b_k^{z|12}=0.
\end{equation}
The characterization multiplicative is given in 
anticipation of the form of the reduced model for $x_1, x_2$ which will have noise terms 
of multiplicative type.

The third model is a combination of the additive and the multiplicative case and 
is called the combined case 
\begin{equation}
\label{com}
\begin{split}
\frac{dx_1}{dt}&=\lambda_a\sum_k b_k^{1|yz}y_k z_k+
\lambda_m\sum_k (b_k^{1|2y}x_2 y_k+b_k^{1|2z}x_2 z_k), \\ 
\frac{dx_2}{dt}&=\lambda_a\sum_k b_k^{2|yz}y_k z_k+
\lambda_m\sum_k (b_k^{2|1y}x_1 y_k+b_k^{2|1z}x_1 z_k), \\
\frac{dy_k}{dt}&=-Re\frac{ik}{2}\sum_{|k'| \leq \Lambda} u_{k'} u_{k-k'}+
\lambda_a b_k^{y|1z}x_1 z_k+\lambda_a b_k^{y|2z}x_2 z_k+
\lambda_m b_k^{y|12}x_1 x_2, \\ 
\frac{dz_k}{dt}&=-Im\frac{ik}{2}\sum_{|k'| \leq \Lambda} u_{k'} u_{k-k'}+
\lambda_a b_k^{z|1y}x_1 y_k+\lambda_a b_k^{z|2y}x_2 y_k+
\lambda_m b_k^{z|12}x_1 x_2,
\end{split}
\end{equation}
where the interaction coefficients are of order 1 and 
satisfy
\begin{equation}
\label{comc}
\begin{split}
b_k^{1|yz}+b_k^{y|1z}+b_k^{z|1y}&=0, \quad b_k^{2|yz}+b_k^{y|2z}+b_k^{z|2y}=0,\\
b_k^{1|2y}+b_k^{2|1y}+b_k^{y|12}&=0, \quad b_k^{1|2z}+b_k^{2|1z}+b_k^{z|12}=0.
\end{split}
\end{equation}
The values of the interaction coefficients used in the numerical experiments can be
found in \cite{majda2}. Due to the constraints on the interaction coefficients and the 
incompressibility of the Burgers bath (when $u_0=0$), the additive model (\ref{add}) 
admits the invariant density 
\begin{equation}
\label{addden}
f(x)=Z^{-1}\exp(-\beta(x_1^2+y_1^2+z_1^2+\ldots+y_{\Lambda}^2+z_{\Lambda}^2)),
\end{equation}
where $Z$ is the normalization constant and $\beta$ is a parameter to be determined 
through the equivalence of the canonical and microcanonical ensembles for large 
enough $\Lambda.$  For the multiplicative and combined cases the density 
\begin{equation}
\label{mulden}
f(x)=Z^{-1}\exp(-\beta(x_1^2+x_2^2+y_1^2+z_1^2+\ldots+y_{\Lambda}^2+z_{\Lambda}^2)),
\end{equation}
is invariant and $Z$ and $\beta$ are defined as in the additive case. To conform with 
the notation of the previous sections we set $\hat{x}=x_1$ in the additive case or 
$\hat{x}=(x_1,x_2)$ in the multiplicative and combined cases, while for all three 
cases $\tilde{x}=(y_1,z_1,\ldots,y_{\Lambda},z_{\Lambda}).$ The form of the density 
implies equipartition of energy among all the variables present in a model (we avoid 
differentiating among the models when there is no risk of confusion) and thus 
$$Var(x_1)=Var(x_2)=Var(y_1)=\ldots=Var(z_{\Lambda})=\frac{1}{2\beta}.$$ For all the numerical 
simulations $\beta=50$ and $\Lambda=50.$ These values are consistent with an initial 
condition with energy $E_0=x_1(0)^2+\sum_{k=1}^{\Lambda}y_1(0)^2+z_1(0)^2$ since 
$$\frac{1}{2\beta}=\frac{E_0}{1+2\Lambda}$$ for the additive model. For the multiplicative 
and combined cases we have $$\frac{1}{2\beta}=\frac{E_0}{2+2\Lambda},$$ where 
$E_0=x_1(0)^2+x_2(0)^2+\sum_{k=1}^{\Lambda}y_1(0)^2+z_1(0)^2.$ In the numerical simulations 
the initial conditions were sampled from the invariant density (\ref{addden}) or 
(\ref{mulden}) (depending on the model at hand) by the Box-Mueller method \cite{kloeden} 
and invoking the equivalence of the canonical and microcanonical ensembles. This 
equivalence was shown for the Burgers bath only in \cite{majda4} for $\Lambda=50$, 
but the form of the invariant density suggests it should hold also for the coupled systems
(see also discussion in \cite{majda1}). This is supported by the fact that, 
in the numerical simulations, the quantities of interest were computed by 
averaging over initial conditions and are identical to the results in 
\cite{majda2}, where they were obtained by time-averaging over a single 
trajectory.

\section{Construction and application of the reduced models}{\label{reduce}}

We proceed to construct the reduced models for the test cases mentioned above.

\subsection{Asymptotic mode reduction strategy models} 
For AMRS the reduced models were constructed in \cite{majda2}, so we only have 
to construct here the reduced models for the short-memory MZ approximation. However, as promised 
in Section \ref{amrs}, we mention briefly the methods involved in the first 
step of AMRS (the stochastic approximation step) for the determination of 
the parameters introduced by the stochastic approximation (\ref{ou}). The parameters introduced 
are the entries of the matrices $\Gamma, \sigma.$ For our test cases, 
$\Gamma=diag(\gamma_1,\ldots,\gamma_{\Lambda})$ and 
$\sigma=diag(\sigma_1,\ldots,\sigma_{\Lambda})$. In the numerical simulations, 
the resolved variables will be coupled only to the first five modes of the Burgers 
bath. However, this does not change anything about the method used to determine 
the parameters of the stochastic approximation. Since we have two parameters, the 
statistical agreement between $$\frac{d\phi}{dt}=R(\phi)$$   
and 
\begin{equation}
\label{odesou}
\begin{split}
\frac{d\hat{\phi}}{dt}&=\hat{R}(\phi) \\
\frac{d\tilde{\phi}}{dt}&=\tilde{R}({\phi})=H(\phi)-
\frac{\Gamma}{\epsilon}\tilde{\phi}+
\frac{\sigma}{\sqrt{\epsilon}} \frac{d\vec{W}(t)}{dt},
\end{split}
\end{equation}
will have to involve the one-time and two-time statistics of the system (see 
\cite{majda2} for details). First, we 
use the one-time statistics to reduce the number of parameters. The unresolved variables evolve on a faster time-scale 
than the resolved ones and, thus, are considered to relax quickly to a stationary state. 
For the stochastic approximation (\ref{ou}), the unresolved variables $y_k, z_k$ reach 
a stationary state  described by a Gaussian density with zero mean and variance 
$\frac{\sigma_k^2}{2\gamma_k}, k=1,\ldots,\Lambda.$ 
On the other hand, we know from (\ref{addden}) or (\ref{mulden}) that every unresolved
variable is Gaussianly distributed with zero mean and variance $\frac{1}{2\beta}.$ Thus,
the identification $\frac{\sigma_k^2}{2\gamma_k}=\frac{1}{2\beta}$ 
yields a relation between $\sigma_k$ and $\gamma_k$ and 
reduces the parameters to be determined to only $\gamma_k$. To do this, we exploit the 
two-time statistics of the system. For the stochastic approximation (\ref{ou}), the 
correlations exhibit exponential decay and the characteristic decay time is $\gamma_k^{-1}.$ 
Of course, the coupling with the resolved variables will alter the decay rates so that the 
correlations do not decay exactly as exponentials. So, the first step is to approximate the 
actual correlation functions by exponentials 
$$\frac{\exp(-\gamma_k^{dns}|t|)}{2\beta},$$ 
where $\gamma_k^{dns}$ is the inverse area 
below the correlations $\mathbf{E}[y_k(t)y_k(0)], \mathbf{E}[z_k(t)z_k(0)]$ normalized 
by $2\beta$:
$$\gamma_k^{dns}=(2\beta \times \text{area under the actual correlation of mode k})^{-1}.$$ 
The superscript denotes that the correlations are computed from simulations of the 
full system and $\mathbf{E}$ denotes expectation with respect to the density for the 
full system. The second step is to pick the values of $\gamma_k$ in (\ref{odesou}) 
so as to optimize consistency with $\exp(-\gamma_k^{dns}|t|)/2\beta.$ 
In \cite{majda2}, there are three different procedures suggested to achieve this. The first 
procedure (P1) uses the scaling law for the correlation times for the Burgers bath \cite{majda4} 
to identify the coefficients $\gamma_k$
$$\gamma_k=\frac{C_1 k}{\sqrt{\beta}},$$
where $C_1$ is a constants to be determined by simulations of the Burgers bath. The 
second procedure (P2) sets $\gamma_k=\gamma_k^{dns}$. The third procedure (P3) adjusts 
the $\gamma_k$ in such a way that the correlations for (\ref{odesou}) reproduce the 
functions $\exp(-\gamma_k^{dns}|t|)/2\beta$ as close as possible (more details 
about the procedures can be found in \cite{majda2}). For the additive case 
the procedure P3 is superior, while for the multiplicative and combined cases the 
procedures P2 and P3 are superior and give similar results. The numerical values 
that we used for the $\gamma_k$ can be found in \cite{majda2}.

After we determine the parameters of the stochastic approximation we can apply the 
elimination procedure outlined in Section \ref{amrs} to obtain the reduced system 
for the resolved variables. The elimination procedure requires the introduction of 
a small parameter $\epsilon$ that controls the magnitude of the different terms that 
appear in the equations. For the additive and multiplicative cases, 
$$\epsilon=\frac{\lambda}{\gamma_{k=1}\sqrt{2\beta}},$$   
while for the combined case
$$\epsilon=\frac{\max(\lambda_a,\lambda_m)}{\gamma_{k=1}\sqrt{2\beta}}.$$ 
The next step is to coarse-grain the time ($t \rightarrow \epsilon t$) and 
apply the second step of the elimination procedure presented in Section \ref{amrs}.

For the 
additive case the reduced equations read
\begin{equation}
\label{addred}
\begin{split}
\frac{dx_1}{dt}&=-\gamma x_1 +\sigma \dot{W}(t)\\
\gamma&=\frac{\lambda^2}{4\beta}\sum_{k=1}^{\Lambda}\frac{(b_k^{1|yz})^2}{\gamma_k},
\quad \sigma=\sqrt{\frac{\gamma}{\beta}}.
\end{split}
\end{equation}
Thus, AMRS predicts that the variable $x_1$ is an Ornstein-Uhlenbeck process.

For the multiplicative case, the reduced equations for the resolved variables read
\begin{equation}
\label{mulred}
\begin{split}
\frac{dx_1}{dt}&=-\lambda^2 \bar{\gamma} x_1-\lambda^2 N_1 x_2^2 x_1 +
\lambda \bar{\sigma}_{11} x_2 \dot{W}_1(t)+\lambda \bar{\sigma}_{12} x_2 \dot{W}_2(t) \\
\frac{dx_2}{dt}&=-\lambda^2 \bar{\gamma} x_2-\lambda^2 N_2 x_1^2 x_2 +
\lambda \bar{\sigma}_{21} x_1 \dot{W}_1(t)+\lambda \bar{\sigma}_{22} x_1 \dot{W}_2(t) \\
\end{split}
\end{equation}
Thus, AMRS predicts that the reduced equations for the resolved variables require 
the introduction of nonlinear terms and multiplicative noises. 
Here $W_1(t),W_2(t)$ are independent Wiener processes and the various parameters 
are defined as follows. Let
\begin{equation}
\label{mulpar1}
\begin{split} 
A&=\beta^{-1}\sum_{k=1}^{\Lambda}\gamma_k^{-1}((b_k^{1|2y})^2+(b_k^{1|2z})^2), \quad
B=\beta^{-1}\sum_{k=1}^{\Lambda}\gamma_k^{-1}((b_k^{2|1y})^2+(b_k^{2|1z})^2) \\
C&=\beta^{-1}\sum_{k=1}^{\Lambda}\gamma_k^{-1}(b_k^{1|2y}b_k^{2|1y}+b_k^{1|2z}b_k^{2|1z}).
\end{split}
\end{equation}
Then
\begin{equation}
\label{mulpar2}
\bar{\gamma}=-\frac{1}{2}C, \quad N_1=\beta(A+C), \quad N_2=\beta(B+C),
\end{equation}
and the matrix $\bar{\sigma}$ is defined as
\begin{equation}
\label{mulpar3}
\bar{\sigma}=\begin{pmatrix} \bar{\sigma}_{11} & \bar{\sigma}_{12} \\
\bar{\sigma}_{21} & \bar{\sigma}_{22} \end{pmatrix},
\end{equation}
and has the property
\begin{equation}
\label{mulpar4}
\bar{\sigma}\bar{\sigma}^T=\begin{pmatrix} A & C \\ C & B \end{pmatrix}.
\end{equation}
The matrix in (\ref{mulpar4}) is positive definite and its square root 
$\bar{\sigma}$ exists, although it is not unique. For the simulations, 
$\bar{\sigma}$ was chosen to be symmetric.

Finally, for the combined case, the reduced equations read
\begin{equation}
\label{comred}
\begin{split}
\frac{dx_1}{dt}=&-\lambda_m^2 \bar{\gamma} x_1-\lambda_m^2 N_1 x_2^2 x_1 +
\lambda_m \bar{\sigma}_{11} x_2 \dot{W}_1(t)+\lambda_m \bar{\sigma}_{12} x_2 \dot{W}_2(t) \\
&-\gamma_{11} x_1 -\gamma_{12} x_2+\sigma_{11} \dot{W}_3(t)+\sigma_{12} \dot{W}_4(t),\\
\frac{dx_2}{dt}=&-\lambda_m^2 \bar{\gamma} x_2-\lambda_m^2 N_2 x_1^2 x_2 +
\lambda_m \bar{\sigma}_{21} x_1 \dot{W}_1(t)+\lambda_m \bar{\sigma}_{22} x_1 \dot{W}_2(t) \\ 
&-\gamma_{12} x_1 -\gamma_{22} x_2+\sigma_{12} \dot{W}_3(t)+\sigma_{22} \dot{W}_4(t).
\end{split}
\end{equation}
The parameters $\bar{\gamma}, N_1,N_2, \bar{\sigma}$ are given by (\ref{mulpar2}-\ref{mulpar4})
and the parameters for the contributions of the additive terms are given by
\begin{equation}
\label{compar2}
\gamma_{11}=\frac{\lambda_a^2}{4\beta}\sum_{k=1}^{\Lambda}\frac{(b_k^{1|yz})^2}{\gamma_k}, \quad
\gamma_{22}=\frac{\lambda_a^2}{4\beta}\sum_{k=1}^{\Lambda}\frac{(b_k^{2|yz})^2}{\gamma_k}, \quad
\gamma_{12}=\frac{\lambda_a^2}{4\beta}\sum_{k=1}^{\Lambda}\frac{b_k^{1|yz}b_k^{2|yz}}{\gamma_k}.
\end{equation}
The matrix $\sigma$ is defined as
\begin{equation}
\label{compar3}
\sigma=\begin{pmatrix} \sigma_{11} & \sigma_{12} \\
\sigma_{21} & \sigma_{22} \end{pmatrix},
\end{equation}
and has the property
\begin{equation}
\label{compar4}
\sigma\sigma^T=\beta^{-1}\begin{pmatrix} \gamma_{11} & \gamma_{12} \\ \gamma_{12} & \gamma_{22} 
\end{pmatrix},
\end{equation}
and like $\bar{\sigma}$ it was chosen to be symmetric for the simulations. Finally, 
the processes 
$W_1(t),W_2(t),W_3(t),W_4(t)$ are independent 
Wiener processes.

\subsection{Short-memory MZ approximation models}

We construct the short-memory MZ equations for the three test cases.  
The invariant probability density used to compute all the necessary projections is 
given by (\ref{addden}) or (\ref{mulden}) depending on the test case. 
The equations for the 
short-memory approximation are 
\begin{equation}
\label{mzs}
\frac{\partial}{\partial{t}} e^{tL}x_j=
e^{tL}PLx_j+e^{tL}QLx_j+
\int_0^t e^{(t-s)L}PLQe^{sL}QLx_jds,
\end{equation}
for $j=1$ or $j=1,2$ depending on the case. For the systems we examine the Markovian 
term $PLx_j$ is identically zero, so 
$QLx_j=Lx_j.$ In this special case, $e^{tQL}QLx_j=e^{tL}Lx_j.$ 
The conditional expectation in the memory term is approximated by 
a finite-rank projection on an orthonormal basis consisting of modified Hermite 
polynomials. The choice of Hermite polynomials was motivated by the Gaussian form
of the density. The basis is given by 
\begin{equation}
\label{her} 
h^{\kappa}(\hat{x})=\prod_{j=1}^{m}\tilde{H}_{\kappa_j},
\end{equation}
where 
\begin{equation}
\label{her2}
\tilde{H}_{\kappa_j}=H_{\kappa_j}((2\beta(1+2\alpha_{\kappa_j}))^{\frac{1}{2}}x_j)
(1+2\alpha_{\kappa_j})^{\frac{1}{4}}
e^{-\alpha_{\kappa_j} \beta x_j^2}.
\end{equation}
For the additive case $m=1$, while for the multiplicative and combined cases $m=2$ and 
the multi-index $\kappa=(\kappa_1,\ldots,\kappa_m).$ The functions $H_{\kappa_j}$ are 
Hermite polynomials (with weight $\exp(-\frac{1}{2}x_j^2))$ satisfying
\begin{equation} 
H_0(x_j)=1, \quad H_1(x_j)=x_j, \quad H_{\kappa_j}(x_j)=\frac{1}{\sqrt{\kappa_j}}
x_jH_{\kappa_j-1}(x_j)-\sqrt{\frac{\kappa_j-1}{\kappa_j}}H_{\kappa_j-2}(x_j)
\label{her3}
\end{equation}
The derivatives of the functions 
$\tilde{H}_{\kappa_j}(x)$ can be computed by the recursive relation
\begin{equation}
\frac{d}{dx}\tilde{H}_{\kappa_j}(x_j)=\sqrt{\kappa_j}
(2\beta(1+2\alpha_{\kappa_j}))^{\frac{1}{2}}
\tilde{H}_{\kappa_{j-1}}(x_j)-2\alpha_{\kappa_j}\beta x_j \tilde{H}_{\kappa_j}(x_j).
\label{her4}
\end{equation}
In general, we allow the different resolved variables to have different scaling 
factors $\alpha_{\kappa_j}$, but in all the numerical simulations the scaling 
factors were set to zero.

We also employ the proposition in 
\cite{stinis} that, if the probability density is invariant  
then the operator $L$ is skew-symmetric. For the case of the finite-rank projection we find
\begin{equation}
\label{mzs2}
\frac{\partial}{\partial{t}} e^{tL}x_j=-\int_0^t\sum_{\kappa \in I} 
(e^{sL}Lx_j,Lh^{\kappa})h^{\kappa}(\hat{\phi}(t-s))ds + e^{tL}Lx_j,
\end{equation}
for $j=1$ or $j=1,2$ depending on the case and $I$ denotes the $m$-tuples of 
indices used in the finite rank projection. Thus, the projection 
coefficients of the memory are correlations of different order of the 
noise term. The projection coefficients in 
the memory term in (\ref{mzs2}) are computed by sampling the invariant density, 
evolving the full system and averaging. The noise term is computed using the 
moving average method for sampling stationary stochastic processes with given 
mean and autocorrelation \cite{gikh1}.

An important issue is how does one decide how many terms and of what form are needed in the expansion of the memory term. For cases, where a scale separation between resolved and unresolved variables exists and the resolved variables include all the slow variables, one way to choose the terms in the expansion is to retain only terms whose coefficients decay 
fast. The larger is the ratio of scale separation in the system  under investigation, 
the more accurate should this approximation be. One would think that maybe what is needed for a better model is the inclusion of slowly decaying memory terms. 
However, this is not necessarily so unless one is willing to solve the orthogonal dynamics equation. The reason is that slowly decaying memory terms 
need to be computed through the solution of the orthogonal dynamics equation (see also Section \ref{smop}). Otherwise, their inclusion can hamper the 
accuracy of the model (evidence of that for the Kuramoto-Sivashinsky equation was presented in \cite{stinis}). In other words, when faced with a 
system that exhibits time scale separation and all the slow variables are resolved, then the appropriate choice for the memory terms (as far as expense/accuracy is 
concerned) is to pick the fast decaying ones. In the limit of infinite scale separation the choice of fast decaying terms becomes exact and this is the 
idea behind the construction of AMRS and short-memory MZ.

The problem with the suggestion of picking only the fast decaying terms is that it quickly becomes impractical. To determine how fast the different coefficients  decay we need to compute them and determine which ones decay fast by visual inspection. For a system with even a few resolved variables the number of possible combinations of basis functions grows very fast with the order of the polynomials used. Additionally, we do not expect all of the fast decaying coefficients to be equally important for the accuracy of the reduced model. Thus, if we include them all, we may reduce the efficiency of the the model unnecessarily.

We propose here a partial fix of the problem of picking the coefficients that are most relevant for systems which exhibit time-scale separation. If the resolved variables include all the slow variables we expect, as already mentioned, the dominant coefficients in the expansion of the memory term to be fast decaying. However, among these fast decaying coefficients the most important ones are those with the largest values at $s=0$ (see (\ref{mzs2})). In other words, if we compute the projection coefficients $(Lx_j,Lh^{\kappa}),$ we suggest to keep in the expansion only the terms whose coefficient magnitude is substantially larger (e.g. at least a couple of orders of magnitude) than the magnitude of the coefficients of the rest of the terms. The advantage of this approach is that we can decide which terms' coefficients we want to compute by a relatively cheap calculation since we do not need to evolve the full system. We applied this criterion to the test cases and concluded that the most important terms in the expansion are similar to the ones predicted by AMRS. The calculation of the evolution of the coefficients for these terms (by integrating repeatedly the full system and averaging) showed that they are indeed fast decaying. However, note that there is the possibility that a projection coefficient $(Lx_j,Lh^{\kappa})$ may have a large value at $s=0$, but the magnitude of the corresponding term $h^{\kappa}$ in the expansion can be small, thus making the above procedure not completely failproof.

We should add here, that the test cases we examine do not exhibit extremely large scale separation ratios and this results in errors in the predictions. This is to be expected and drives home the point that for an improved prediction we should, also, include in the expansion slowly-decaying coefficients of higher order, computed through the orthogonal dynamics equation. In MZ, in addition to picking up only the fast decaying terms in the expansion, we have at our disposal the orthogonal dynamics  equation whose solution, although expensive, can provide us with slowly decaying projection coefficients even for the case where there is no scale separation between the 
resolved and unresolved variables. Of course, this ability stems from the fact that MZ is a reformulation of the problem, hence it allows, in principle, the construction of models of arbitrary accuracy. In AMRS, it is not possible, except for special cases (see Section 4.5 in \cite{majda1}), to account for slowly decaying terms.

It is true that the MZ formalism allows refinements by computing the
evolution of slowly decaying memory terms. This can be done through
the orthogonal dynamics equation. However, there are two sources of
additional computational expense in this case. One is to solve the orthogonal
dynamics equation and the other is to solve the
random integrodifferential equation for the reduced model, which
will now have more integral terms. The solution of the orthogonal dynamics equation is not particularly
expensive because the quantities needed to set it up
can be computed in the same process of computing the quantities needed
for the short-memory model. Of course, there is an additional expense to
actually solve the orthogonal dynamics equation but this is not as severe. The
reason is that the solution of the orthogonal dynamics equation can be
reduced to a solution of a system of linear Volterra integral equations which
can be done efficiently (see also \cite{chorin1}). The main expense appears after one has constructed the reduced model.
First, there are more integral terms now than before, and moreover, these terms
have slowly decaying integrands. As a result, the repeated solution
of such a reduced model can become prohibitively expensive because the
evaluation of the integral terms can be very costly.

\subsection{Numerical simulations}{\label{numbers}}

We present the results of the numerical simulation of the reduced models 
produced by AMRS and MZ for the test cases of Section \ref{models}. We should note here that the amount of work needed to construct and integrate the reduced models is comparable for the 
two methods. For AMRS, one needs to solve the full system repeatedly to compute the entries of
the dissipation matrix $\Gamma$. For short-memory MZ, one needs to solve the full system repeatedly to 
compute the memory kernels $(e^{sL}Lx_j,Lh^{\kappa}).$ After the reduced models are 
computed, one needs to integrate them. In the case of AMRS, we need to integrate a system
of stochastic differential equations. In the case of MZ, a system of random integrodifferential 
equations. The fact that the MZ reduced equations are random (colored noise) allows the faster convergence of the averaging procedure (e.g. to compute correlations). However, the gain in efficiency is offset by the fact that we deal with integrodifferential equations. As a result, the overall computational effort needed is the same for the two methods. More details about the implementation of the reduced models are offered below when the examples are analyzed.

\begin{figure}
\centering
\subfigure[]{\epsfig{file=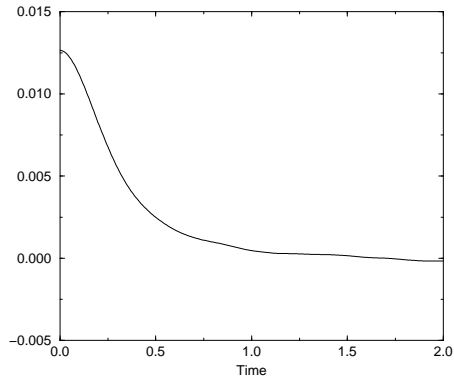,height=2in}}
\caption{Additive case. Evolution of $(e^{tL}Lx_1, L\tilde{H}_1(x_1)).$}
\label{fig:add1}
\end{figure}

For the additive case the AMRS reduced equation for the resolved variable $x_1$ is
given by (\ref{addred}) and the MZ reduced equation is given by
\begin{equation}
\label{addopred}
\begin{split}
\frac{d\phi_1}{dt}&=-\int_0^{t_0} (e^{sL}Lx_1,L\tilde{H}_1(x_1))\tilde{H}_1(\phi_1(t-s))ds
+F_1(t) \\
\phi_1(0)&=x_1(0),
\end{split}
\end{equation}
where 
$$L=R_{x_1}\pd{}{x_1}+\sum_{k=1}^{\Lambda}(R_{y_k}\pd{}{y_k}+R_{z_k}\pd{}{z_k}), 
\quad Lx_1=R_{x_1}=\lambda\sum_k b_k^{1|yz}y_k z_k,$$ 
and $F_1(t)$ is a stationary stochastic process with mean zero and autocorrelation 
$(e^{tL}Lx_1, Lx_1),$ where the inner product is weighted by the invariant 
density (\ref{addden}). The coupling constant $\lambda=4,$ the truncation size 
is $\Lambda=50$ and $\beta=50.$ Note that due to 
the short-memory approximation, the 
interval of the integral in (\ref{addopred}) is restricted to $[0,t_0].$ This is
done for two reasons. First, because the kernel of the integral $(e^{sL}Lx_1,L\tilde{H}_1(x_1))$ 
decays fast and second, because the estimate of the kernel can not be accurate 
for large $s$ by the error estimates for the short-memory approximation. The 
temporal evolution of the kernel is shown in Fig.(\ref{fig:add1}). It was 
computed by sampling the density (\ref{addden}), evolving the full system (\ref{add}) 
with the fourth-order Runge-Kutta algorithm and averaging. The results shown in 
Fig.(\ref{fig:add1}) correspond to averaging over 10000 samples. For the numerical 
simulations of (\ref{addopred}) we set $t_0=1.$ The fact that the kernel decays fast
motivates its replacement by a delta-function multiplied by the 
integral of the kernel. This approximation will be referred to as the delta MZ 
approximation. To check the relevance of the short-memory approximation, we also 
solved the orthogonal dynamics equation for $x_1$ keeping up to second order terms 
in the expansion of the memory and the temporal evolution of 
$(e^{tQL}Lx_1, L\tilde{H}_1(x_1))$ is identical to 
$(e^{tL}Lx_1, L\tilde{H}_1(x_1)).$ 
The parameter 
$\epsilon=0.63$ for this case, which is not small, but as seen from the results 
below it is small enough to guarantee the validity of the short-memory approximation 
and of AMRS.

\begin{figure}
\centering
\subfigure[]{\epsfig{file=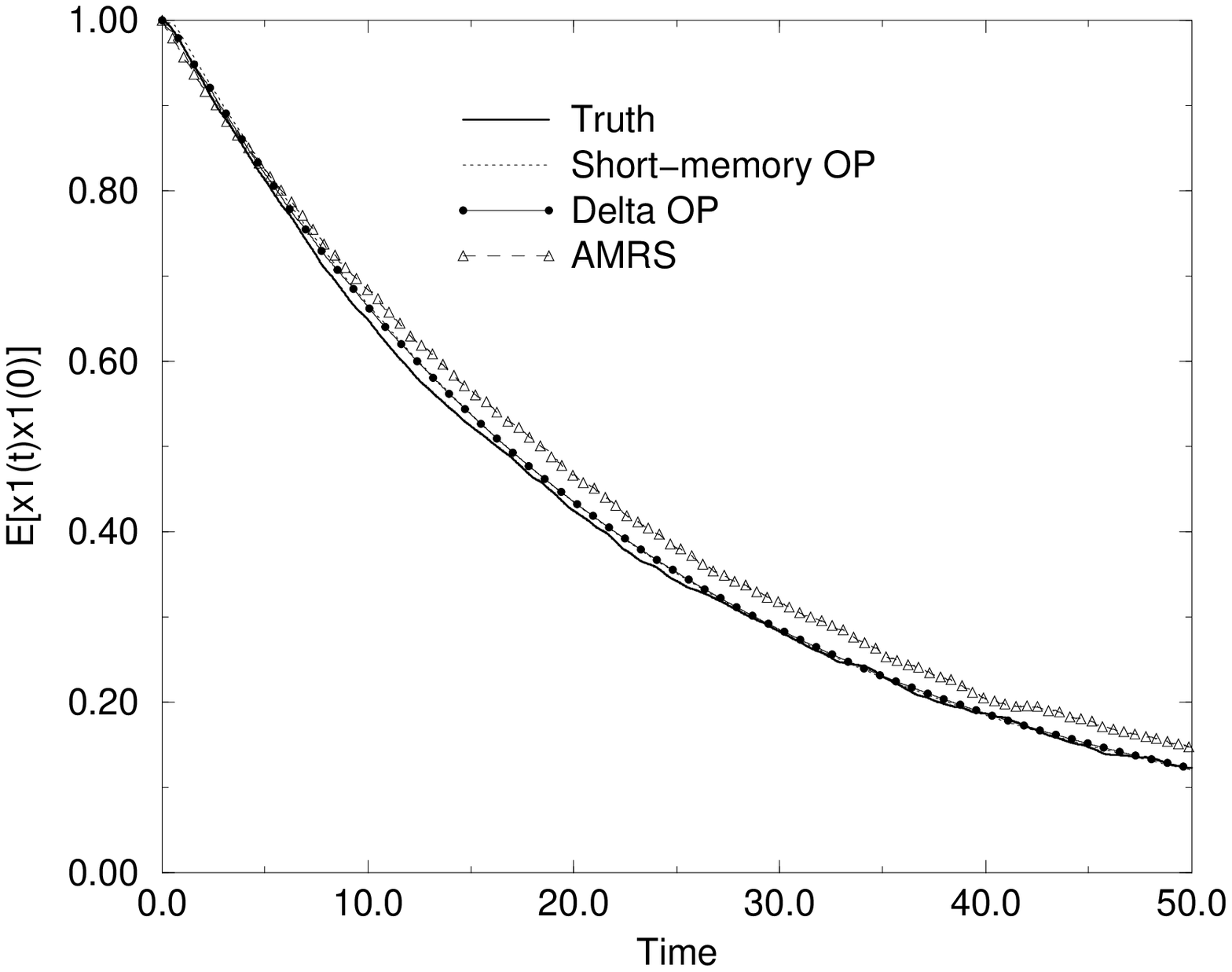,height=2in}}
\quad
\subfigure[]{\epsfig{file=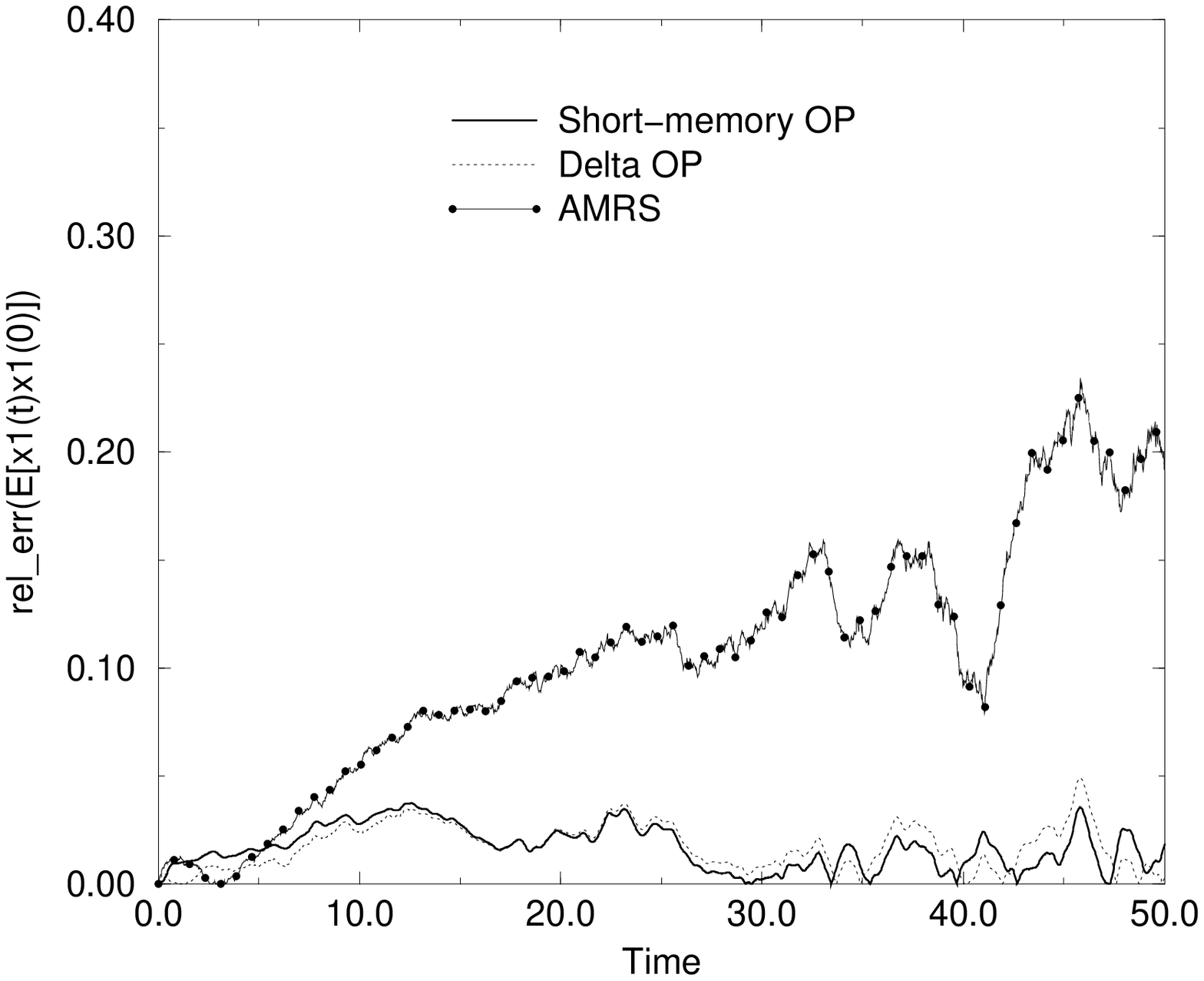,height=2in}}
\caption{Additive case. a) Autocorrelation of the resolved mode $x_1.$ 
b) Relative error of the predictions of the autocorrelation of $x_1.$ }
\label{fig:add2}
\end{figure}

In \cite{majda2}, the autocorrelation of $x_1$ was used to evaluate the performance of the 
reduction strategy. Fig.(\ref{fig:add2}) shows the predictions for the autocorrelation 
of $x_1$ from AMRS, short-memory MZ and the delta MZ approximation. The truth refers to 
the autocorrelation computed from the full system (in this case (\ref{add})). 
As mentioned before, 
only the first five interactions coefficients were nonzero. The truth was computed 
by sampling the density (\ref{addden}), evolving (\ref{add}) with the fourth-order 
Runge-Kutta algorithm and averaging. The estimate was computed using 10000 samples. 
The AMRS reduced equation (\ref{addred}) was evolved with the Euler scheme for initial 
values of $x_1$ 
drawn from the projected on $x_1$ form of the density (\ref{addden}) and different 
realizations of the white noise term. The AMRS estimate 
of the autocorrelation was computed by averaging over 20000 samples and noise realizations. 
Finally, the MZ estimates were produced by sampling again the projected on $x_1$ density 
and the noise term $F_1(t)$, 
evolving with the fourth-order Runge-Kutta algorithm and averaging. Note that, the MZ
equations have a colored noise term and thus do not require the use of a stochastic 
solver. The short-memory MZ estimates were computed by averaging over 1000 samples and noise 
realizations. The use of a colored noise allows for the faster convergence of the 
averaging procedure. However, the gain in computational time by the need to 
evolve fewer samples is balanced by the fact that we have to solve integrodifferential 
equations for the short-memory approximation. The delta MZ approximation, where 
we no longer have integral memory terms, performs surprisingly well, since it has almost 
the same accuracy as the short-memory MZ approximation while being much more efficient 
numerically (about 10 times). The delta MZ approximation was computed by averaging over 
1000 samples and noise realizations. All three estimates have comparable accuracy for 
short times, while for later times the relative error of the AMRS estimate appears to 
increase reaching a plateau of about $20\%.$ The MZ estimates' relative error remains 
around $3\%$ for the interval of integration.

\begin{figure}
\centering
\subfigure[]{\epsfig{file=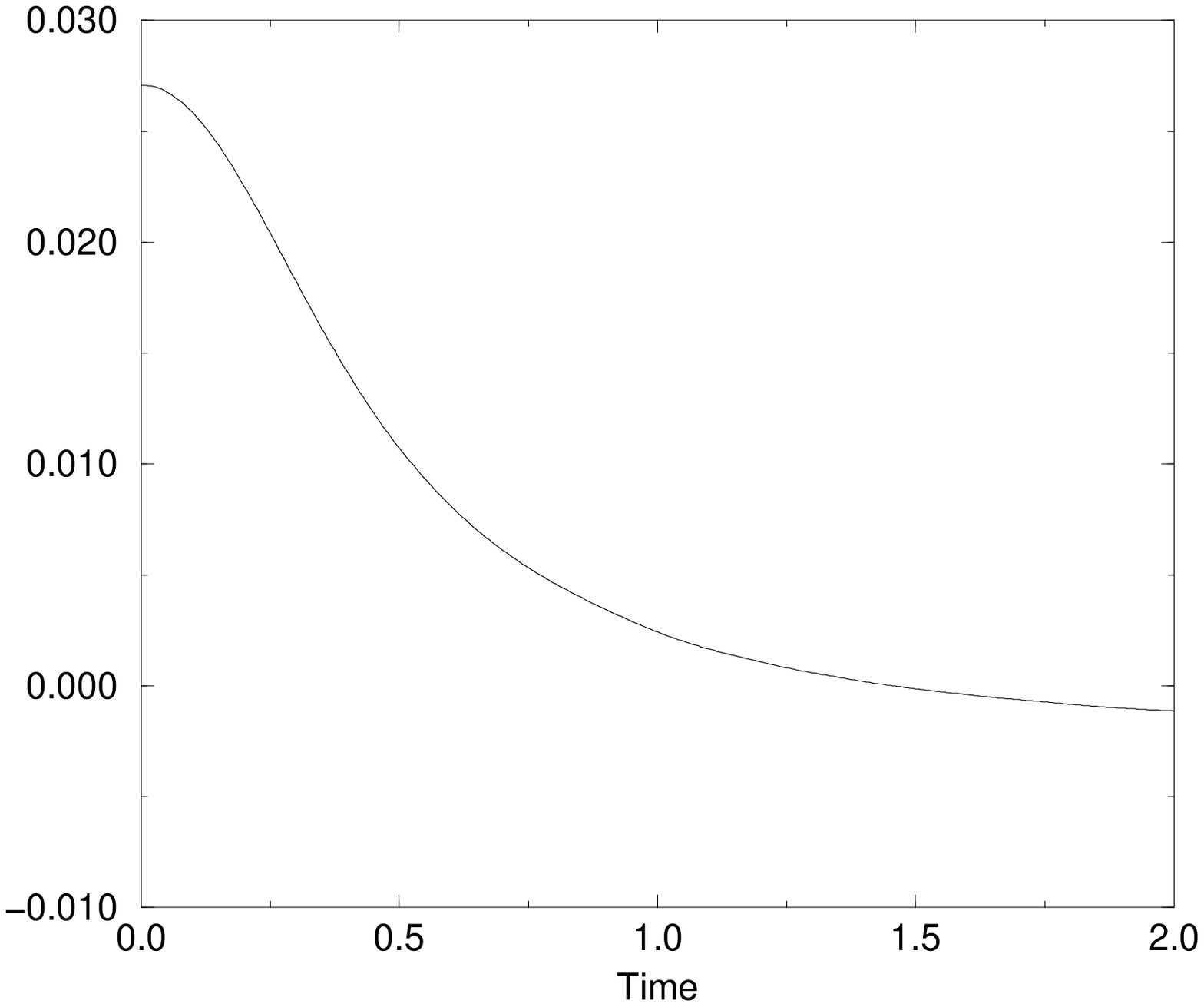,height=2in}}
\quad
\subfigure[]{\epsfig{file=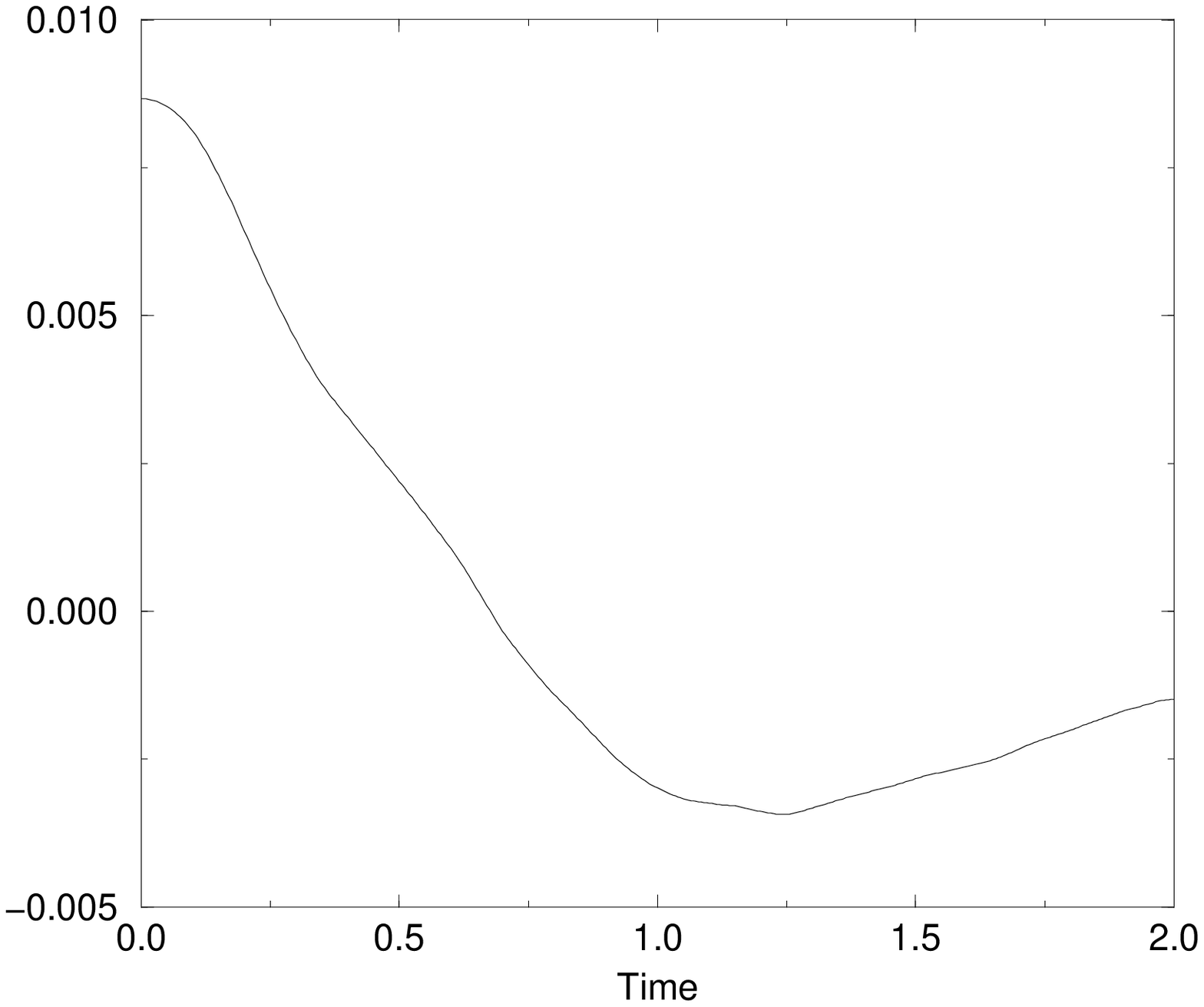,height=2in}}
\caption{ Multiplicative case. a) Evolution of $(e^{tL}Lx_1, L\tilde{H}_1(x_1)).$
b) Evolution of $(e^{tL}Lx_1,L(\tilde{H}_2(x_2)\tilde{H}_1(x_1))).$}
\label{fig:mul1}
\end{figure}

For the multiplicative case (see (\ref{mul})), the AMRS reduced equations for $x_1,x_2$ 
are given by (\ref{mulred}) while the MZ equations are 
\begin{equation}
\label{mulopred}
\begin{split}
\frac{d\phi_1}{dt}=&-\int_0^{t_0} (e^{sL}Lx_1,\tilde{H}_1(x_1))\tilde{H}_1(\phi_1(t-s))ds\\
&-\int_0^{t_0}(e^{sL}Lx_1,L(\tilde{H}_2(x_2)\tilde{H}_1 (x_1)))
\tilde{H}_2(\phi_2(t-s))\tilde{H}_1(\phi_1(t-s))ds+F_1(t), 
\\
\frac{d\phi_2}{dt}=&-\int_0^{t_0} (e^{sL}Lx_2,L\tilde{H}_1(x_2))\tilde{H}_1(\phi_2(t-s))ds \\
&-\int_0^{t_0}(e^{sL}Lx_2,L(\tilde{H}_2(x_1)\tilde{H}_1(x_2)))
\tilde{H}_2(\phi_1(t-s))\tilde{H}_1(\phi_2(t-s))ds+F_2(t), 
\\
\phi_1(0)=&x_1(0), \quad \phi_2(0)=x_2(0),
\end{split}
\end{equation}
where 
\begin{equation*}
\begin{split}
L&=R_{x_1}\pd{}{x_1}+R_{x_2}\pd{}{x_2}+
\sum_{k=1}^{\Lambda}(R_{y_k}\pd{}{y_k}+R_{z_k}\pd{}{z_k}), \\ 
Lx_1&=R_{x_1}=\lambda\sum_k (b_k^{1|2y}x_2 y_k+b_k^{1|2z}x_2 z_k), \\ 
Lx_2&=R_{x_2}=\lambda\sum_k (b_k^{2|1y}x_1 y_k+b_k^{2|1z}x_1 z_k),
\end{split}
\end{equation*} 
and $F_1(t),F_2(t)$ are stationary stochastic processes with mean zero and autocorrelation 
$(e^{tL}Lx_1, Lx_1),
(e^{tL}Lx_2, Lx_2)$ respectively. The inner product is weighted by 
the invariant density (\ref{mulden}). The 
coupling constant $\lambda=3.$ The interval 
of integration is again restricted to $[0,t_0]$ and for the numerical simulations $t_0=2.$ 
Fig.(\ref{fig:mul1}) shows the temporal evolution of the memory kernels 
$(e^{tL}Lx_1, L\tilde{H}_1(x_1))$ and  
$(e^{tL}Lx_1,L(\tilde{H}_2(x_2)\tilde{H}_1(x_1))).$ The kernels for the equation 
for $x_2$ have similar behaviour. 
The kernels decay fast and, as in the previous case, we also tried to replace the 
kernels by a delta function multiplied by the integral of the kernel and this approximation 
is again called delta MZ. However, the fact that we have to integrate the kernels up to 
$t_0=2$ shows that they decay more slowly than in the additive case. This is to be expected
if we look at the form of, say, $Lx_1.$ Due to the multiplicative coupling, this 
quantity now depends not only on the unresolved variables, but, also, on the slow resolved 
variable $x_2$ and thus the autocorrelation of $Lx_1$ depends on the autocorrelation of 
$x_2$. This results in a slower decay of $(e^{tL}Lx_1, L\tilde{H}_1(x_1))$, and we expect 
the error of the short-memory approximation to be larger than in the additive case. On 
the other hand, AMRS results in equations for the resolved variables that have a 
multiplicative noise (see (\ref{mulred})), i.e. the dependence on the slow variables is taken out of the noise process. 
This is in accordance with the way the AMRS equations are derived, which is in the limit of 
infinite-scale separation.

The parameter $\epsilon=0.41.$ All algorithmic 
considerations and numbers 
of samples are the same as in the additive case, except for the algorithm used to 
integrate the AMRS reduced equations (\ref{mulred}). We followed \cite{majda2} and 
used time-splitting, using a second-order Runge-Kutta algorithm for the nonlinear 
terms and the strong Milstein scheme of order 1 for the stochastic terms.

\begin{figure}
\centering
\subfigure[]{\epsfig{file=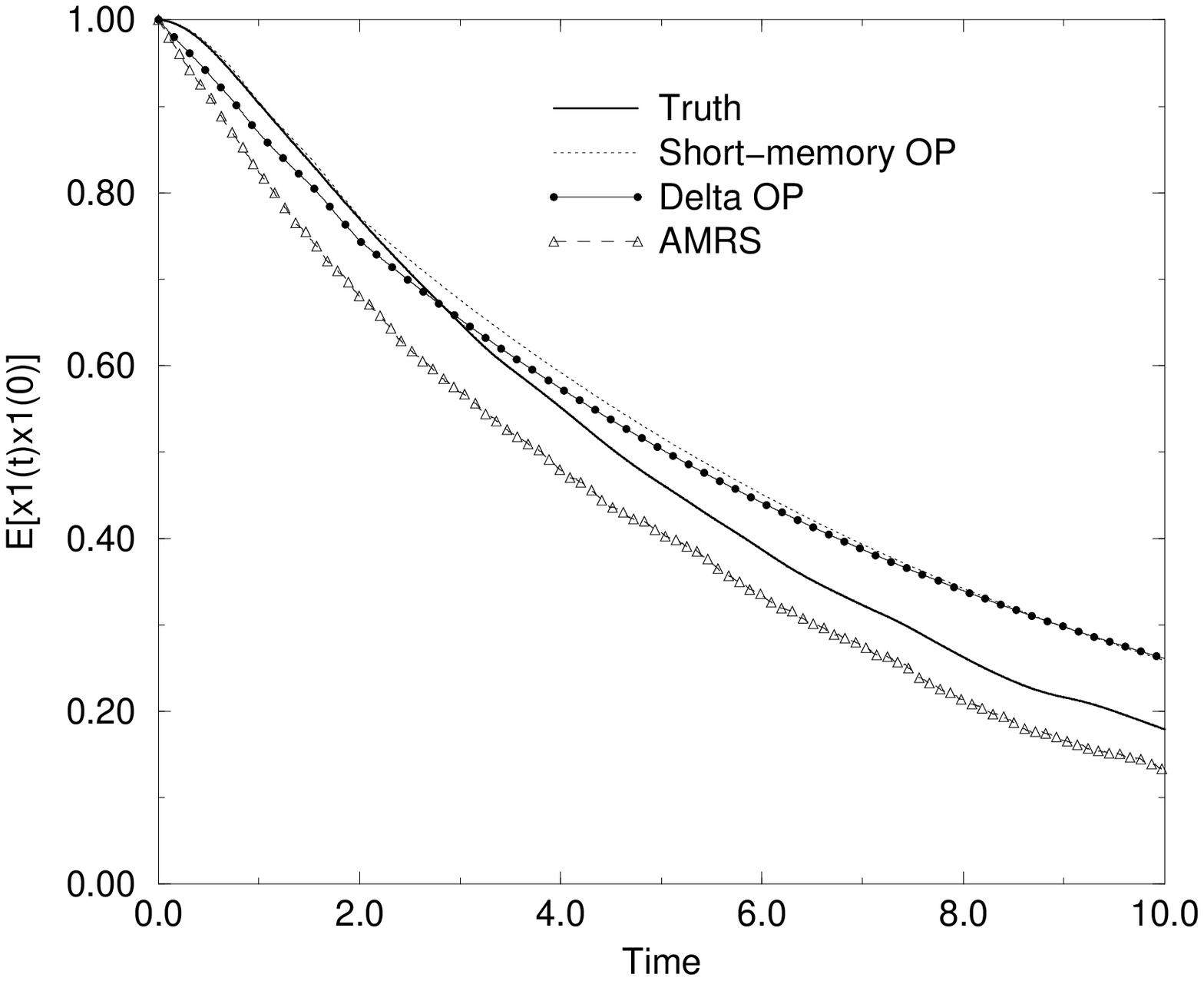,height=2in}}
\quad
\subfigure[]{\epsfig{file=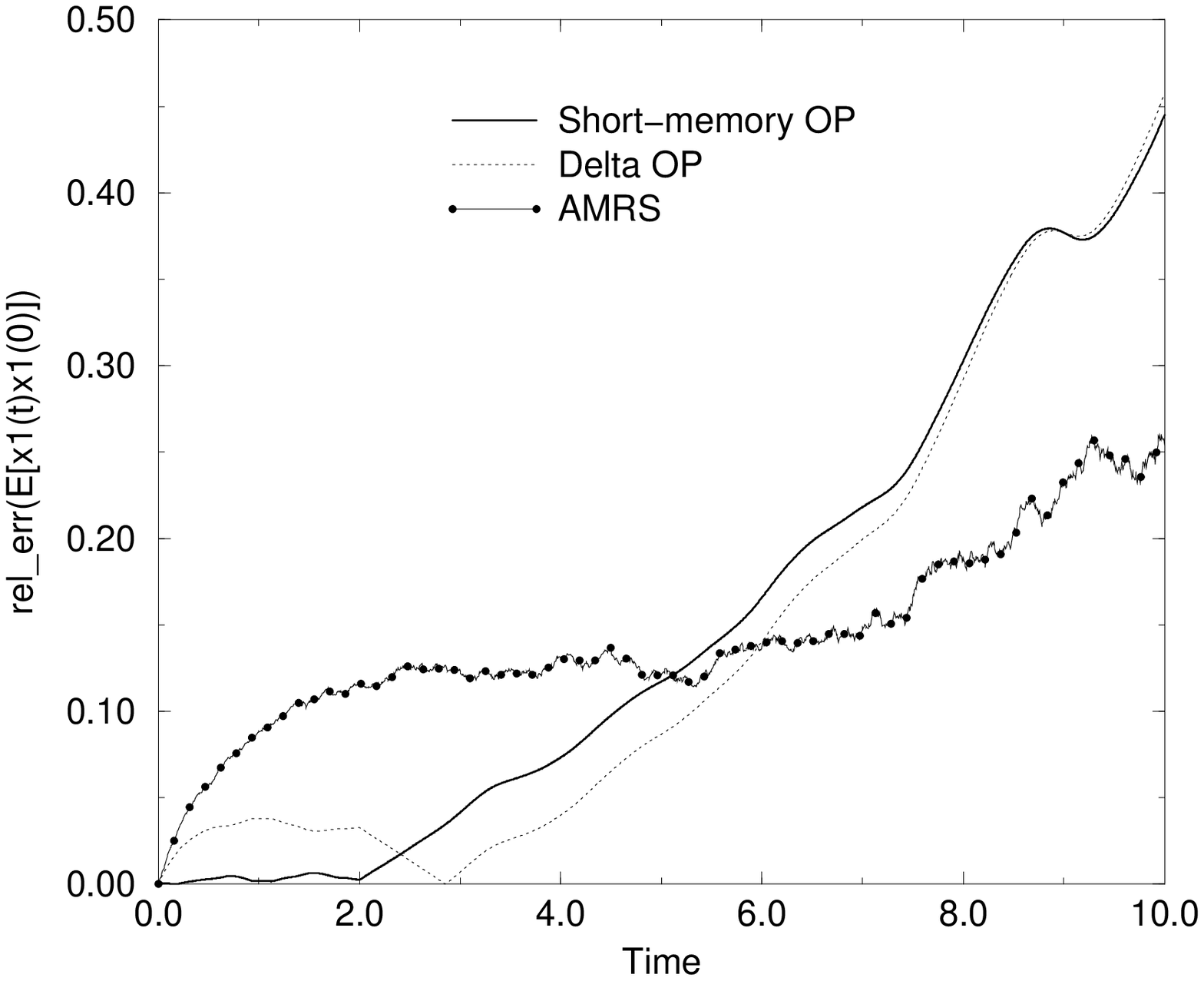,height=2in}}
\caption{Multiplicative case. a) Autocorrelation of the resolved mode $x_1.$ 
b) Relative error of the predictions of the autocorrelation of $x_1.$ }
\label{fig:mul2}
\end{figure}

To check again the validity of the short-memory approximation we solved the orthogonal 
dynamics equations for $x_1,x_2$ keeping up to second order polynomials in each 
variable and the results for the kernels were very close to the ones 
predicted by the short-memory approximation. This suggests that we expect the short-memory and AMRS estimates to be, at least for short times, accurate. Fig.(\ref{fig:mul2}) shows
the estimates for the autocorrelation of $x_1$ as predicted by AMRS, short-memory MZ 
and delta MZ compared to the truth (similar results hold for $x_2$). The MZ estimates 
are superior to the ones produced by AMRS for up to time $t\approx 5,$ after which the 
AMRS estimate becomes more accurate. Note that, up to time $t\approx 5$ the error for 
all methods is around $10\%,$ thus justifying our expectation of the validity 
of the short-memory approximation for short times. At the instant $t=10$ when we stopped 
the integration the relative error of the AMRS estimate is about $25\%$ while the one 
produced by the two MZ variants is about $45\%$. The larger error compared to the additive 
case confirms our expectations stated above. The fact that delta MZ and short-memory MZ 
give very similar results and that to within second order polynomials the orthogonal dynamics 
kernels are close to the short-memory ones, indicates that the short-memory kernels are rather 
accurate. So, what seems to be needed for an improvement of the results is to 
compute, through the orthogonal dynamics, projection coefficients of higher order 
and, possibly, for these coefficients 
extend the interval of time integration in the memory term. 
We did not attempt this here, because our 
purpose was to compare AMRS and MZ in situations where their numerical efficiency is 
of the same order.

\begin{figure}
\centering
\subfigure[]{\epsfig{file=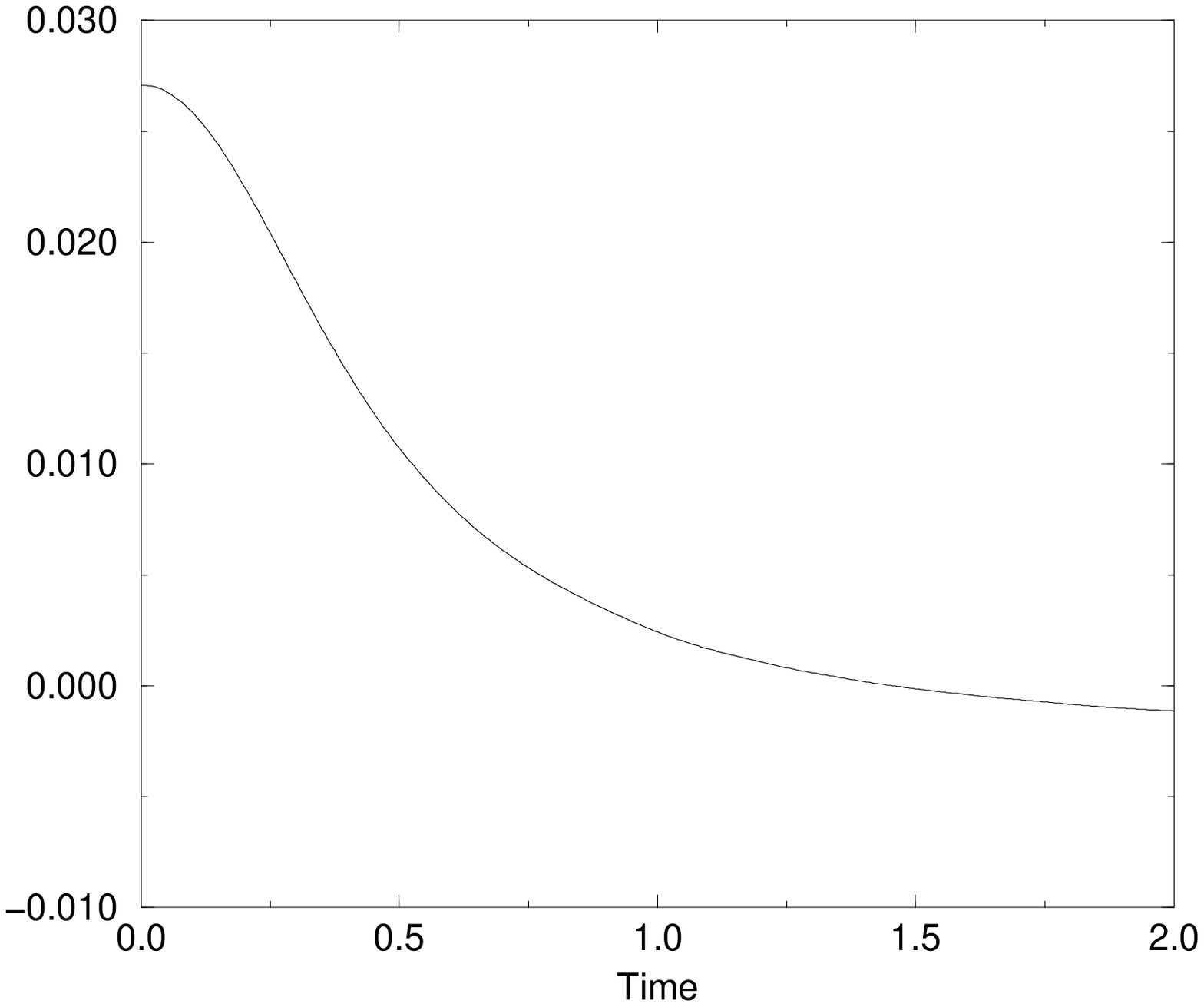,height=2in}}
\quad
\subfigure[]{\epsfig{file=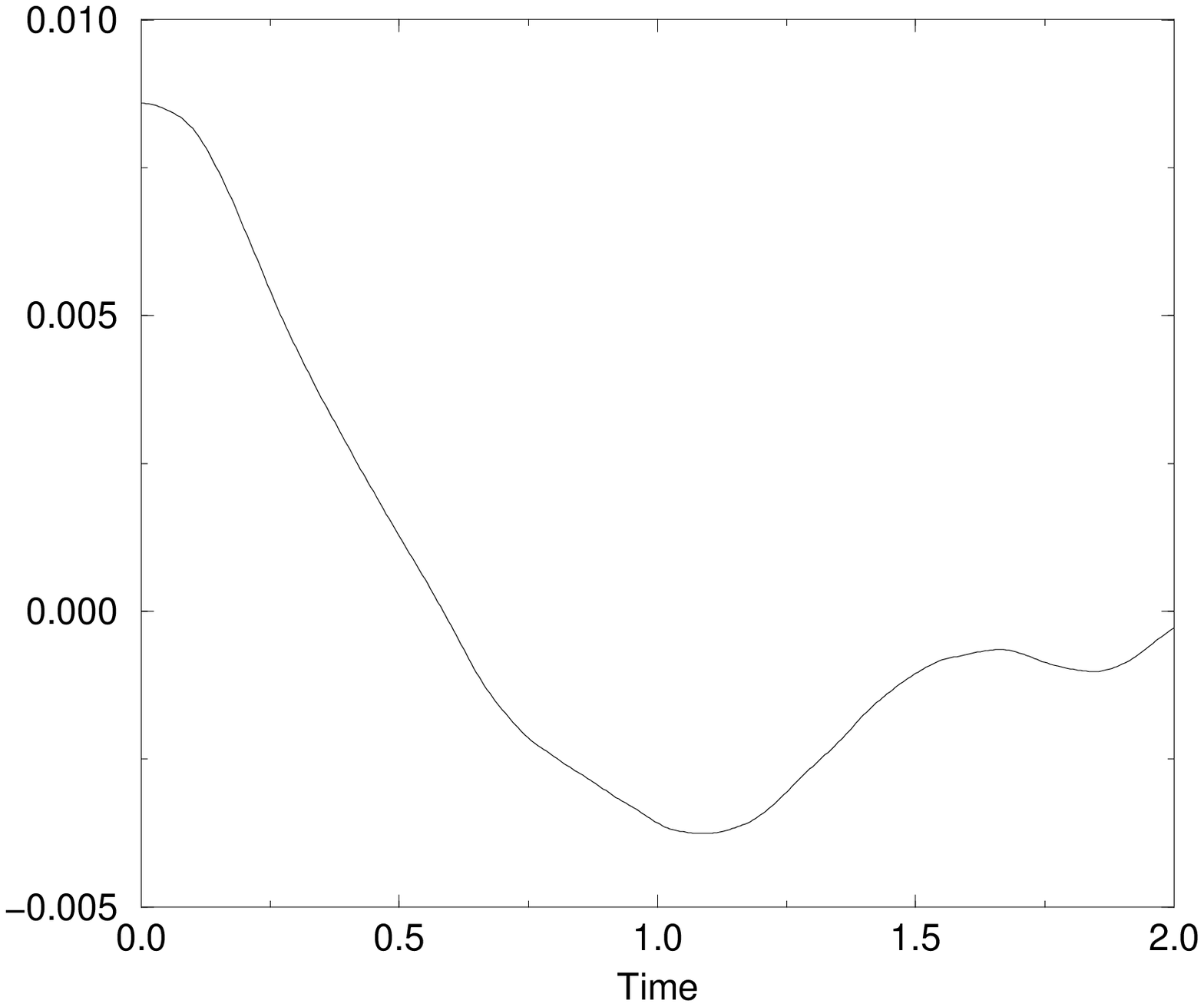,height=2in}}
\caption{ Combined case. a) Evolution of $(e^{tL}Lx_1, L\tilde{H}_1(x_1)).$
b) Evolution of $(e^{tL}Lx_1,L(\tilde{H}_2(x_2)\tilde{H}_1(x_1))).$}
\label{fig:com1}
\end{figure}

The last test case we examined is the combined one (see (\ref{com})), where the coupling of 
the resolved variables to the unresolved is of both additive and multiplicative type. 
For this case, the AMRS reduced equations are given by (\ref{comred}), while the 
reduced MZ equations are 
\begin{equation}
\label{comopred}
\begin{split}
\frac{d\phi_1}{dt}=&-\int_0^{t_0} (e^{sL}Lx_1,\tilde{H}_1(x_1))\tilde{H}_1(\phi_1(t-s))ds\\
&-\int_0^{t_0}(e^{sL}Lx_1,L(\tilde{H}_2(x_2)\tilde{H}_1 (x_1)))
\tilde{H}_2(\phi_2(t-s))\tilde{H}_1(\phi_1(t-s))ds+F_1(t), 
\\
\frac{d\phi_2}{dt}=&-\int_0^{t_0} (e^{sL}Lx_2,L\tilde{H}_1(x_2))\tilde{H}_1(\phi_2(t-s))ds \\
&-\int_0^{t_0}(e^{sL}Lx_2,L(\tilde{H}_2(x_1)\tilde{H}_1(x_2)))
\tilde{H}_2(\phi_1(t-s))\tilde{H}_1(\phi_2(t-s))ds+F_2(t), 
\\
\phi_1(0)=&x_1(0), \quad \phi_2(0)=x_2(0),
\end{split}
\end{equation}
where 
\begin{equation*}
\begin{split}
L&=R_{x_1}\pd{}{x_1}+R_{x_2}\pd{}{x_2}+
\sum_{k=1}^{\Lambda}(R_{y_k}\pd{}{y_k}+R_{z_k}\pd{}{z_k}), \\ 
Lx_1&=R_{x_1}=\lambda_a\sum_k b_k^{1|yz}y_k z_k+
\lambda_m\sum_k (b_k^{1|2y}x_2 y_k+b_k^{1|2z}x_2 z_k), \\ 
Lx_2&=R_{x_2}=\lambda_a\sum_k b_k^{2|yz}y_k z_k+
\lambda_m\sum_k (b_k^{2|1y}x_1 y_k+b_k^{2|1z}x_1 z_k),
\end{split}
\end{equation*} 
and $F_1(t),F_2(t)$ are stationary stochastic processes with mean zero and autocorrelation 
$(e^{tL}Lx_1, Lx_1),(e^{tL}Lx_2, Lx_2)$ respectively. The inner product is weighted by 
the invariant density (\ref{mulden}). The 
coupling constants are $\lambda_a=4, \lambda_m=3.$ The interval 
of integration is again restricted to $[0,t_0]$ and for the numerical simulations $t_0=1.$ 
Note, that while the reduced AMRS equations (\ref{comred}) involve terms that 
come either from the multiplicative coupling part or the additive coupling part, the 
MZ equations (\ref{comopred}) involve also cross-terms, i.e. products 
of terms where one factor comes from the multiplicative part and the other from 
the additive part. Fig.(\ref{fig:com1}) shows the temporal evolution of the kernels 
$(e^{tL}Lx_1, L\tilde{H}_1(x_1)),(e^{tL}Lx_1,L(\tilde{H}_2(x_2)\tilde{H}_1(x_1))).$ 
The kernels for the equation for 
$x_2$ have similar behaviour. 
The kernels decay fast and, as in the previous case, we also tried to replace the 
kernels by a delta function multiplied by the integral of the kernel and this approximation 
is again called delta MZ. The parameter $\epsilon=0.49.$ All algorithmic 
considerations and numbers 
of samples are the same as in the multiplicative case.

\begin{figure}
\centering
\subfigure[]{\epsfig{file=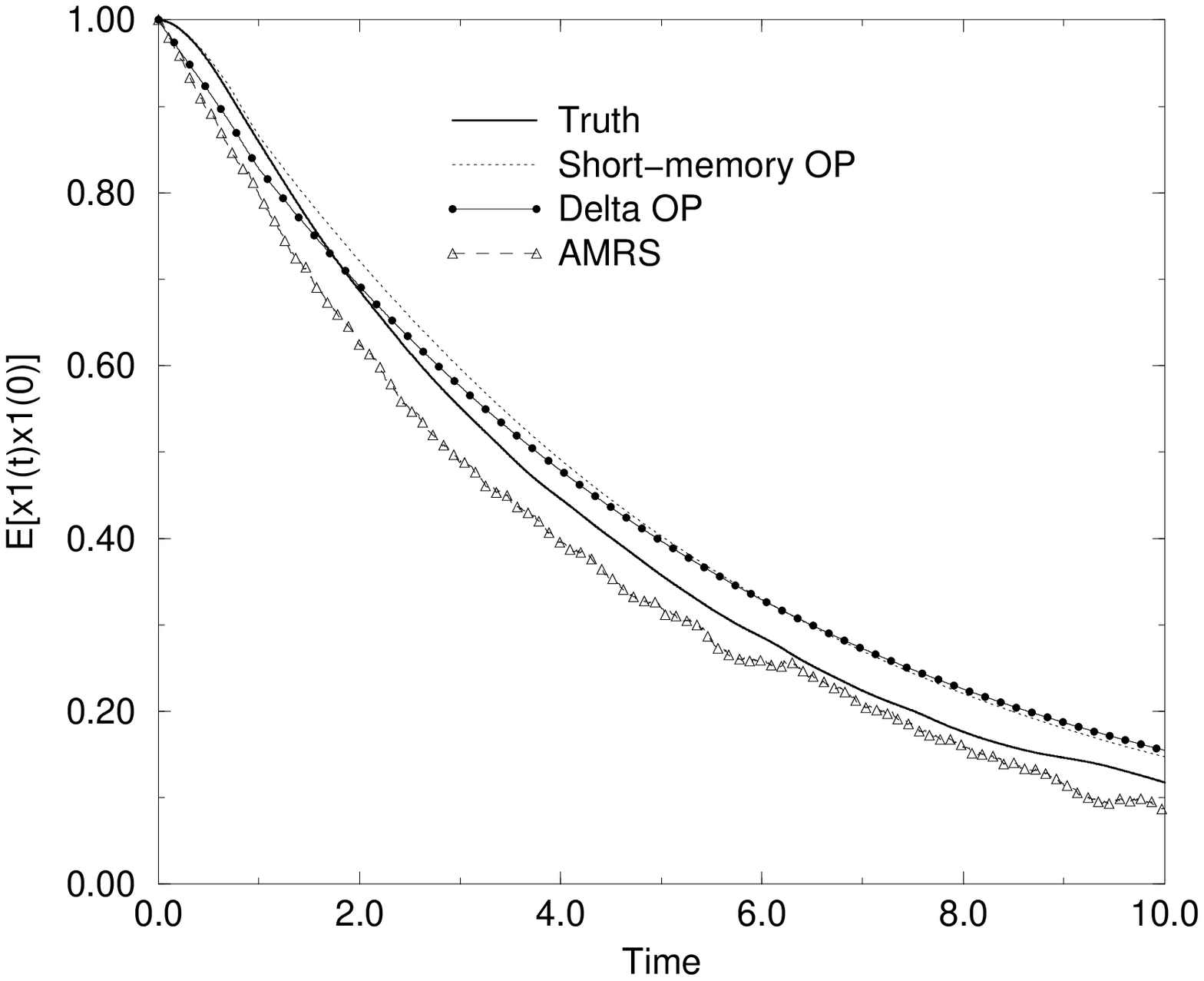,height=2in}}
\quad
\subfigure[]{\epsfig{file=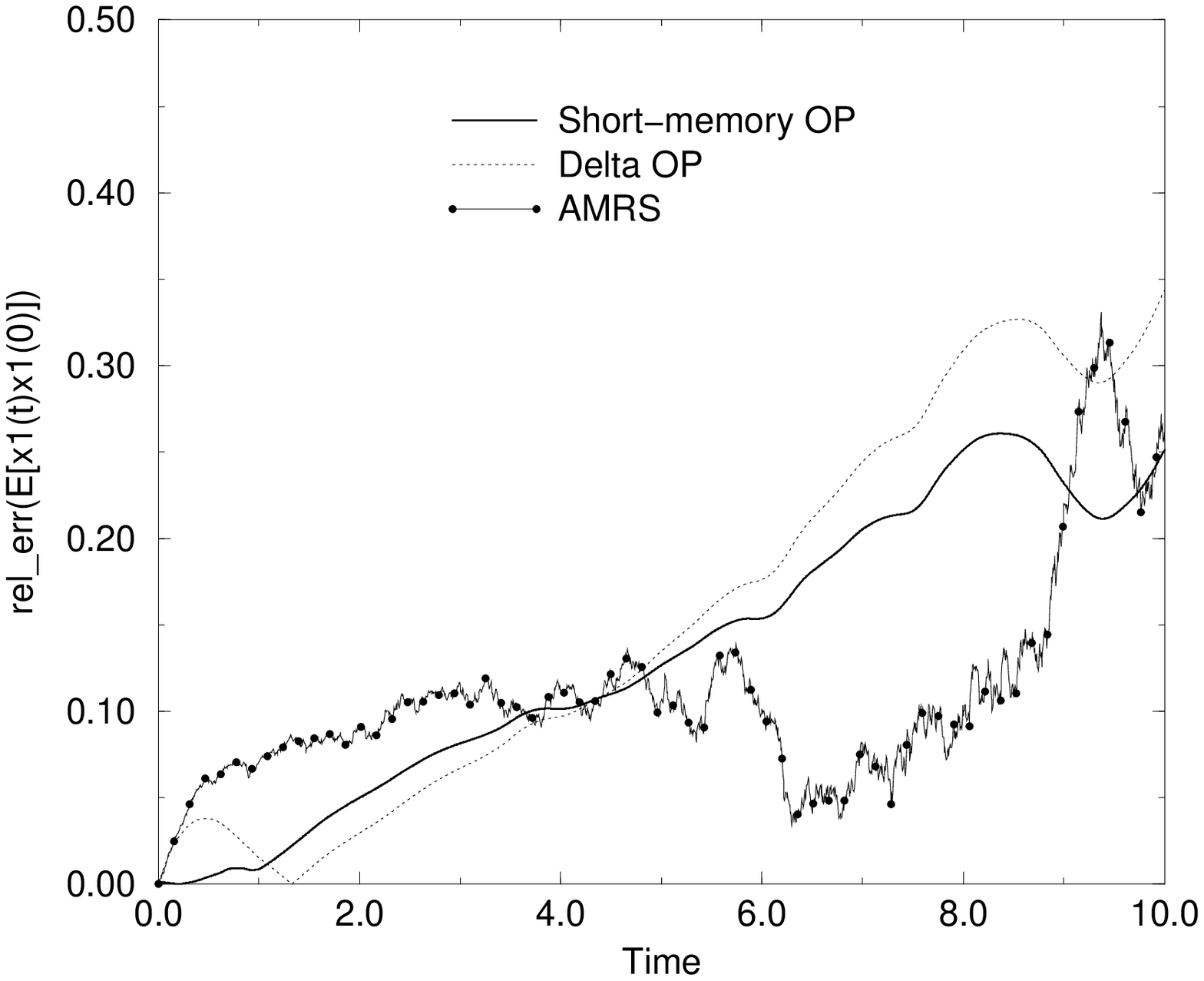,height=2in}}
\caption{Combined case. a) Autocorrelation of the resolved mode $x_1.$ 
b) Relative error of the predictions of the autocorrelation of $x_1.$ }
\label{fig:com2}
\end{figure}

Fig.(\ref{fig:com2}) shows
the estimates for the autocorrelation of $x_1$ as predicted by AMRS, short-memory MZ 
and delta MZ compared to the truth (similar results hold for $x_2$). The MZ estimates 
are superior to the ones produced by AMRS for up to time $t\approx 5.$ After that and 
until time $t \approx 9,$ the AMRS estimate's accuracy is higher. For the time 
interval shown here, the relative error for all three methods has a maximum 
value around $25\%.$ If we compare the 
accuracy of the results for MZ and AMRS taking into account their behaviour in the 
pure additive and multiplicative cases, we see that the performance in the combined 
case is to be expected, since it is a mixture of the advantages and drawbacks of 
each individual method as manifested in the additive and multiplicative cases. The 
fact, that MZ exhibits higher accuracy for a larger interval of time, should 
be attributed to the fact that the additive coupling constant is larger, thus bringing 
the combined case somewhat closer to the additive case (see also discussion of results 
for the combined case in \cite{majda2}).

\section{Conclusions}{\label{conc}}
We have presented numerical results comparing two stochastic mode reduction strategies. 
The first method (AMRS), proposed by Majda, Timofeyev and Vanden-Eijnden,
is based on an asymptotic strategy developed by Kurtz. 
The second method is a short-memory approximation of the Mori-Zwanzig projection formalism 
of irreversible statistical mechanics, as proposed by 
Chorin, Hald and Kupferman. The novel feature 
of these methods is that they allow the use, in the reduced system, of higher order 
terms in the resolved variables. The two methods were applied to a collection 
of test cases that exhibit separation of time-scales between the resolved and unresolved 
variables and, also, share some features with more complicated models used in the 
study of climate dynamics.

Depending on the test case, one of the two methods can have 
superior accuracy, but the overall behaviour suggests that for cases with separation 
of time-scales, the two methods result in reduced systems of equations that have similar 
predictive ability. For the test cases we examined, the replacement of the 
kernels in the memory terms by delta-functions (the delta MZ 
approximation) does not appear to be very harmful to the accuracy of the approximation, 
while at the same time it makes the 
integration of the reduced equations around 10 times faster. 
The test cases highlight the 
limitations of AMRS and of short-memory MZ when the separation of time-scales becomes 
less sharp. In this case, MZ allows for a systematic, although expensive, 
calculation of reduced equations that incorporate long-time memory effects. On the other 
hand, AMRS, by construction, cannot be readily extended to these cases. For special 
cases, like the multiplicative case above (see also Section 4.5 in \cite{majda1}), 
the reduction performed by AMRS can be effected by working directly on the 
stochastic differential equations (\ref{rndodes}), without recourse to the 
associated Chapman-Kolmogorov equation. This allows for a systematic development of reduced 
equations that account for long-time memory effects.

It would be interesting to compare the two methods when applied to more realistic 
models e.g. equations for climate dynamics, where a separation of time-scales is 
known to exist between the quantities of interest and the huge number of faster 
variables that constitute the climate system \cite{majda5}.

\section*{Acknowledgements}
I would like to thank Professor Alexandre Chorin and Professor Andrew Majda for very helpful discussions and comments. I would also like to thank the reviewers for their insightful comments. This work was supported in part by the Applied Mathematical Sciences subprogram of the Office of Energy Research of the US Department of Energy under Contract DE-AC03-76-SF00098 and in part by the National Science Foundation under Grant DMS98-14631.

\bibliographystyle{elsart-num}

\bibliography{paper}

\end{document}